\documentclass{amsart}

\usepackage[mathscr]{eucal}
\usepackage{amsmath}
\usepackage{amssymb}
\usepackage{amscd}
\usepackage{amsfonts}
\usepackage{enumerate}

\usepackage[all]{xy} 
\CompileMatrices

\theoremstyle{plain}
\newtheorem{thm}{Donotwrite}[section]

\newtheorem{theorem}[thm]{Theorem}
\newtheorem{proposition}[thm]{Proposition}
\newtheorem{lemma}[thm]{Lemma}
\newtheorem{corollary}[thm]{Corollary}

\newtheorem{algo}[thm]{Rule}

\theoremstyle{definition}
\newtheorem{example}[thm]{Example}

\newtheorem{remark}[thm]{Remark}

\numberwithin{equation}{section}

\newfont{\germ}{eufm10}

\newcommand\La{\Lambda}

\newcommand\lw[1]{\lower.4mm\hbox{${#1}$}}

\newcommand\veps{\varepsilon}
\newcommand\vphi{\varphi}

\newcommand{\be}{\begin{enumerate}}
\newcommand{\ee}{\end{enumerate}}

\newcommand{\eq}{\begin{eqnarray}}
\newcommand{\eneq}{\end{eqnarray}}
\newcommand{\Lemma}{\begin{lemma}}
\newcommand{\enlemma}{\end{lemma}}
\newcommand{\ba}{\begin{array}}
\newcommand{\ea}{\end{array}}

\newcommand{\batten}[5]{%
\begin{picture}(40,40)(-20,-20)
\put(-10,0){\line(1,0){20}}
\thicklines
\put(0,10){\line(0,-1){20}}
\put(-11,0){\makebox(0,0)[r]{$#1$}}
\put(0,16){\makebox(0,0)[b]{$#2$}}
\put(0,-16){\makebox(0,0)[t]{$#3$}}
\put(11,0){\makebox(0,0)[l]{$#4$}}
\end{picture}
}

\title[Scattering Rule in Soliton Cellular Automaton 
associated with Crystal Base of $U_q(D_4^{(3)}) $]
{Scattering Rule in Soliton Cellular Automaton 
associated with Crystal Base of $U_q(D_4^{(3)}) $}

\author{Daisuke Yamada}
\address{Daisuke Yamada: 
Graduate School of Mathematical Sciences, the University of Tokyo, 
Megro, Tokyo 153-8914, Japan} 

\date{}

\begin{document}

\mathversion{bold}
\maketitle
\mathversion{normal}


\begin{abstract}\noindent
In terms of the crystal base of a quantum affine algebra $U_q(\mathfrak{g})$, 
we study a soliton cellular automaton (SCA) 
associated with the exceptional affine Lie algebra $\mathfrak{g}=D_4^{(3)}$. 
The solitons therein are labeled by the crystals of quantum affine algebra 
$U_q(A_1^{(1)})$. 
The scatteing rule is identified with the combinatorial $R$ matrix for 
$U_q(A_1^{(1)})$-crystals. 
Remarkably, the phase shifts in our SCA are given by 
{\em $3$-times} of those in the well-known box-ball system.
\end{abstract}

\section{Introduction}

\subsection{Box ball systems in terms of crystals}

The box-ball system [19, 21] is a discrete dynamical system in which finitely many balls move along the one-dimensional array of boxes under a certain rule. 
Then each isolated array composed of balls behaves like a soliton: 
{\em (i) The longer isolated soliton moves faster.
(ii) The number of solitons does not change under the time evolution. 
(iii) If the solitons are apart enough from each other at an initial state, then their lengths do not change after their collisions.}  
Its integrability has been proved by making a connection to the difference analog of the Lotka-Volterra equation through the limiting procedure called ultra-discretization [22]. 
In the last part of 90's, the box-ball system was re-formulated [8] in terms of crystal base theory [9, 14, 15] associated with the symmetric tensor representations of the quantum affine algebra $U_q(\mathfrak{g})$ for non-exceptional affine Lie algebra $\mathfrak{g}_n=B_n^{(1)}$, $C_n^{(1)}$, $D_n^{(1)}$, $A_{2n-1}^{(2)}$, $A_{2n}^{(2)}$ and $D_{n+1}^{(2)}$. 
Let us recall the box-ball system in terms of crystals. 
For simplicity we restrict the discussion in $\mathfrak{g}_n=A_n^{(1)}$ 
[2, 3]. 
By a similar discussion, 
there are some expansions [7, 24].


Recall necessary results on box ball system in terms of crystals. 
For $U'_q(A_n^{(1)})$-crystal $B_l$ ($l\in \mathbb{Z}_{\geq 1}$), 
see section \ref{sec:An crystal}. 
Fix a sufficiently large positive integer $L$. 
Put $u_l=(l,0,\ldots, 0)\in B_l$.  
Define a set of paths $\mathscr{P}_L$ by 
\begin{equation*}
\mathscr{P}_L=\left\{ \ 
p=b_1\otimes b_2\otimes \cdots \otimes b_L \in B_1^{\otimes L} 
;\ b_j=u_1 \ \text{if $j\gg 1$} \ \right\}. 
\end{equation*} 
{} By the iterating crystal isomorphism $B_l\otimes B_1\simeq B_1\otimes B_l$, 
we have a map from $B_l\otimes \mathscr{P}_L$ to $\mathscr{P}_L\otimes B_l$ by 
\begin{align*}
u_l\otimes b_1 \otimes & b_2\otimes \cdots \otimes b_L 
\\ \mapsto & \ 
\tilde{b}_1 \otimes u'_l\otimes b_2 \otimes  \cdots \otimes b_L
\mapsto \ \cdots \quad \cdots \ \mapsto
\tilde{b}_1 \otimes \tilde{b}_2\otimes \cdots \otimes 
\tilde{b}_L \otimes u_l.
\end{align*} 
{} Thus we define a map $T_l : \mathscr{P}_L\longrightarrow \mathscr{P}_L$ by 
\begin{equation*}
T_l( b_1 \otimes b_2 \otimes \cdots \otimes b_L )=
\tilde{b}_1 \otimes \tilde{b}_2 \otimes \cdots \otimes \tilde{b}_L. 
\end{equation*}
Application of $T_l$ induces a transition of state. 
Thus map $T_l$ can be regarded as a dynamical system, 
in which $T_l$ plays the role of {\em time evolution}. 
We call the system $\{T_l^t(p) \ \vert \ t\geq 0,\ l\in \mathbb{Z}_{\geq 1},\ p\in \mathscr{P}_L \}$ the {\em $A_n^{(1)}$-automaton}. 
Thanks to the Yang-Baxter equation, 
$\{T_l\}_{l\geq 1}$ becomes a commuting family. 

A $m$-soliton state of length $(l_1,l_2,\ldots ,l_m)$ is given by 
\begin{equation*}
...[l_1]......[l_2].....\cdots .....[l_m]..............
\end{equation*}
Here $\ldots [l]\ldots$ denotes a local configuration such as 
\begin{equation}\label{eq:affine A_n-soliton}
\cdots \otimes (1)\otimes (1)
\otimes \underbrace{(n+1)^{\otimes x_{n}}\otimes \cdots \otimes (3)^{\otimes x_2}\otimes (2)^{\otimes x_1}}_{x_1+x_2+\cdots +x_n=l}\otimes (1)\otimes (1)\otimes \cdots 
\end{equation}
sandwiched by sufficiently many $(1)$'s. 
Here $(\cdot )$ signifies a frame of Young tableau. 

The state $T_k(p)$ is obtained by the rightward shift by 
$\min(k,l)$ lattice steps. 
The following figure explains the action $T_l$ ($l=3$) on a 
one-soliton state (It transfers the upper line into the lower line): 
\[
T_3: 3\otimes 3\otimes 2\otimes 1\otimes 1\otimes 1\otimes 1\otimes \cdots
\mapsto 1\otimes 1\otimes 1\otimes 3\otimes 3\otimes 2\otimes 1\otimes \cdots 
\]
\[
\batten{111}{3}{1}{}{a}
\hspace{-1mm}
\batten{113}{3}{1}{}{a}
\hspace{-1mm}
\batten{133}{2}{1}{}{a}
\hspace{-1mm}
\batten{233}{1}{3}{}{a}
\hspace{-1mm}
\batten{123}{1}{3}{}{a}
\hspace{-1mm}
\batten{112}{1}{2}{}{a}
\hspace{-1mm}
\batten{111}{1}{1}{\cdots}{a}
\]
The array ``$332$"  behaves as a soliton just like 
in the classical soliton theory, where the frames are omitted. 

Let $\tilde{B}_l$ ($l\in \mathbb{Z}_{\geq 1}$) be a $U'_q(A_{n-1}^{(1)})$-crystal. Define a map $\iota_l: \tilde{B}_l\sqcup \{0\}\rightarrow B_1^{\otimes l}\sqcup \{0\}$ by 
\begin{equation}
\iota_l(b)=
\begin{cases}
(n+1)^{\otimes x_n}\otimes \cdots \otimes (3)^{\otimes x_2}\otimes (2)^{\otimes x_1} & \text{if $b\neq 0$}, 
\\ 
0 & \text{if $b=0$}. 
\end{cases}
\end{equation}
{} Then it is easy to see that for any $i\in \{1,2,\ldots ,n-1\}$ and 
$b\in \tilde{B}_l\sqcup \{0\}$ 
\begin{equation}
\iota_l(\tilde{e}_i(b))=\tilde{e}_{i+1}(\iota_l(b)),\quad  
\iota_l(\tilde{f}_i(b))=\tilde{f}_{i+1}(\iota_l(b))
\end{equation}
{} Thus local states $\ldots [l]\ldots$ are labeled by crystals of the smaller algebra $U_q(A_{n-1}^{(1)})$. 
We also denote $m$-soliton state with length $(l_1,l_2,...,l_m)$ by 
the tensor product of $U_q(A_{n-1}^{(1)})$-crystals. 
\begin{equation}
z^{\gamma_1}b_1\otimes z^{\gamma_2}b_2\otimes \cdots 
\otimes z^{\gamma_m}b_m \in
\text{Aff}(\tilde{B}_{l_1})\otimes 
\text{Aff}(\tilde{B}_{l_2})\otimes \cdots \otimes \text{Aff}(\tilde{B}_{l_m})
\end{equation}
Here $\gamma$ denotes a phase of the local configuration $...[l]...$ 
and defined by $\gamma=\min(r,l)t+
\text{(the position of $``...[l]..."$ from the right end)}$ under $T_r$. 
It is invariant unless $...[l]...$ interacts with other solitons. 


For the $m$-soliton state with length $(l_1,l_2,\ldots ,l_m)$, 
we assume $l_1>l_2>\cdots >l_m$. 
Since the longer solitons move faster under the time evolution $T_r$ $(r\geq l_1)$, we can expect that the state turns out to be 
\begin{equation}
z^{\gamma'_m}b'_m \otimes \cdots \otimes z^{\gamma'_2}b'_2 \otimes 
z^{\gamma'_1}b'_1\in \text{Aff}(\tilde{B}_{l_m})\otimes \cdots \otimes \text{Aff}(\tilde{B}_{l_2})\otimes \text{Aff}(\tilde{B}_{l_1}).  
\end{equation}
after sufficiently many time evolutions. 
We represent such a scattering process as follows. 
\begin{equation}
S_m: 
z^{\gamma_1}b_1\otimes z^{\gamma_2}b_2 \otimes \cdots \otimes z^{\gamma_m}b_m
\mapsto z^{\gamma'_m}b'_m \otimes \cdots \otimes z^{\gamma'_2}b'_2 \otimes 
z^{\gamma'_1}b'_1
\end{equation}

\begin{theorem}
$([2,\ 3])$\ 
For $m>2$ the map $S_m$ is factorized into $S_2$.
The map $S_2$ is described by the combinatorial $R$ 
for $U_q(A_{n-1}^{(1)})$-crystals. Especially, the phase shift is given by 
\begin{equation}
\gamma'_2-\gamma_2=\gamma_1-\gamma'_1=2l_2+H(b_1\otimes b_2). 
\end{equation}
\end{theorem}

\subsection{Our case}

In [16] a family of $U'_q(D_4^{(3)})$-crystals $\{B_l\}_{l\geq 1}$ corresponding to the vertex $1$ in the Dynkin diagram was studied. 
Kashiwara, Misra, Okado and the author gave an explicit description of the crystal structure by the {\em coordinate representation}  
(see. section \ref{sec:D_4^(3)}) in the spirit of [11]. 
Thus it is expected that a soliton cellular automaton associated with 
$U'_q(D_4^{(3)})$-crystals can be constructed by the same method as 
[3, 8]. 

First, we have to make a table of a combinatorial $R$ matrix 
for $U'_q(D_4^{(3)})$-crystals. 
\begin{eqnarray*}
R: \text{Aff}(B_l)\otimes \text{Aff}(B_1) & \longrightarrow & 
\text{Aff}(B_1)\otimes \text{Aff}(B_l)
\\ 
z^{\gamma_1}b_1\otimes z^{\gamma_2}b_2 & \mapsto & 
z^{\gamma_2+H(b_1\otimes b_2)}b'_2\otimes z^{\gamma_1-H(b_1\otimes b_2)}b'_1. 
\end{eqnarray*}
Now we introduce an algorithm for the combinatorial $R$ (section \ref{sec:comb.R}) by the Lecouvey's column insertion algorithm for type $G_2$ [17]. 
We describe it only for those tableau with shape $(n)$ or $(n+1,1)$ for $n\in \mathbb{Z}_{\geq 0}$, which suffices for our aim. 
The proof goes as follows.   
As $U_q(G_2)$-crystals, $B_l\otimes B_1$ decomposes into connected components that are isomorphic to the crystals of irreducible $U_q(G_2)$-modules. 
Within each component, the general elements are obtained by applying 
$\tilde{f}_i$'s $(i=1,2)$ to the $U_q(G_2)$ highest weight elements.  
Thus we first prove the theorem directly for the highest elements 
(See Appendix A). 
For the general elements, the theorem follows from the fact due to [17]. 
(See Proposition \ref{prop:G2 morphism} and Lemma \ref{lem:G2}.) 

Next, using the algorithm of combinatorial $R$ given in the previous section, 
we construct a soliton cellular automaton associated with 
$U'_q(D_4^{(3)})$-crystals.  We call it the {\em $D_4^{(3)}$-automaton}. 
The method is same as the $A_n^{(1)}$ case. 
We summarize the result below: 
(We use the same notation as the $A_n^{(1)}$ case.)

\begin{itemize}
\item A one-soliton state in $D_4^{(3)}$-automaton is described as 
\begin{equation}
\cdots \otimes (1)\otimes (1)
\otimes \underbrace{(3)^{\otimes x_2}\otimes (2)^{\otimes x_1}}_{x_1+x_2=l}\otimes (1)\otimes  \cdots 
\iff  z^{\gamma}(x_1,x_2)\in \text{Aff}(\tilde{B}_l).
\end{equation}
Here $\tilde{B}_l$ denotes the $U'_q(A_1^{(1)})$-crystal. 

\item The two-body scattering of solitons is described by the combinatorial $R$ for $U'_q(A_1^{(1)})$-crystals. Especially, the phase shift is given by 
\begin{equation}
\gamma'_2-\gamma_2=\gamma_1-\gamma'_1=2l_2+3\times H(b_1\otimes b_2). 
\end{equation}

\item Scattering of solitons is factorized into two-body scattering. 
\end{itemize}

Let us compare the $A_2^{(1)}$-automaton with $D_4^{(3)}$-automaton. 
Let $\{\tilde{B}_l\}_{l\geq 1}$ be the family of $U'_q(A_1^{(1)})$-crystals. 
Under the time evolution $T_2$, 
we give an example of two-soliton scattering described by 
$\text{Aff}(\tilde{B}_2)\otimes \text{Aff}(\tilde{B}_1)\rightarrow 
\text{Aff}(\tilde{B}_1)\otimes \text{Aff}(\tilde{B}_2)$. 
\[
\begin{array}{rl}
{} & \text{$D_4^{(3)}$-automaton}
\\ & \\ 
t=0:& \underline{22}111\underline{3}11111111111111111111
\\ 
t=1:& 11\underline{22}11\underline{3}1111111111111111111
\\ 
t=2:& 1111\underline{22}1\underline{3}111111111111111111
\\ 
t=3:& 111111\underline{22}\underline{3}11111111111111111
\\ 
t=4:& 11111111\underline{2\underline{0}}1111111111111111
\\ 
t=5:& 1111111111\underline{\underline{\bar{3}}1}11111111111111
\\ 
t=6:& 11111111111\underline{\phi}\underline{21}111111111111
\\ 
t=7:& 111111111111\underline{1}0\underline{21}1111111111
\\ 
t=8:& 1111111111111\underline{1}23\underline{21}11111111
\\ 
t=9:& 11111111111111\underline{1}213\underline{21}111111
\\ 
t=10:& 111111111111111\underline{1}2113\underline{21}1111
\\ 
t=11:& 1111111111111111\underline{1}21113\underline{21}11
\end{array}
\begin{array}{l}
{}  \text{$A_2^{(1)}$-automaton}
\\  \\ 
\underline{22}111\underline{3}111111111111111111111
\\ 
11\underline{22}11\underline{3}11111111111111111111
\\ 
1111\underline{22}1\underline{3}1111111111111111111
\\ 
111111\underline{22}\underline{3}111111111111111111
\\ 
11111111\underline{232}1111111111111111
\\ 
1111111112\underline{13}211111111111111
\\ 
11111111112\underline{1}\underline{13}2111111111111
\\ 
111111111112\underline{1}1\underline{13}21111111111
\\ 
1111111111112\underline{1}11\underline{13}211111111
\\ 
11111111111112\underline{1}111\underline{13}2111111
\\ 
111111111111112\underline{1}1111\underline{13}21111
\\ 
1111111111111112\underline{1}11111\underline{13}211
\end{array}
\]
Here the markers specify the positions of initial solitons under the interaction-free propagation. We can read off the phase shift from the deviation of the outgoing solitons from the corresponding markers. In terms of the soliton labels, the above scattering is described by the combinatorial $R$ for $U_q(A_1^{(1)})$-crystal: 
\begin{eqnarray*}
\text{Aff}(\tilde{B}_2)\otimes \text{Aff}(\tilde{B}_1)
& \longrightarrow & 
\text{Aff}(\tilde{B}_1)\otimes \text{Aff}(\tilde{B}_2)
\\ 
z^{L-2}(2,0)\otimes z^{L-6}(0,1) & \mapsto & 
z^{(L-6)+\delta}(1,0)\otimes z^{(L-2)-\delta}(1,1).   
\end{eqnarray*}
Here the phase shift $\delta$ is given by 
\begin{align*}
\delta=&2\times 1+H((2,0)\otimes (0,1))\times 
\begin{cases}
3& \text{if $\mathfrak{g}=D_4^{(3)}$},\\ 
1& \text{if $\mathfrak{g}=A_2^{(1)}$}
\end{cases}
\\=&
-1 \quad \text{if $\mathfrak{g}=D_4^{(3)}$} ;\ =1 \quad \text{if $\mathfrak{g}=A_2^{(1)}$}. 
\end{align*}

In order to prove the theorem, we introduce an operator $T_{\natural}$ 
constructed from the $U'_q(D_4^{(3)})$-crystal $B_{\natural}$ 
which corresponds to the vertex $2$ in the Dynkin diagram. 
The idea is same as [7]. 
$T_{\natural}$ is to put the highest element 
$\begin{pmatrix} 1\\ 2 \end{pmatrix}\in B_{\natural}$ 
on the left and to carry it through to the right via the isomorphism 
$B_{\natural}\otimes B_1^{\otimes L}\simeq B_1^{\otimes L}\otimes B_{\natural}$. It acts on the space of states and has the following properties: 
(See Example \ref{ex:ironuki}.)
\begin{align}
& 
\text{$T_{\natural}$ commutes with $T_r$.} 
\\ & 
\text{$
T_{\natural}(z^{\gamma}(x_1,x_2))=
\begin{cases}
z^{\gamma}(x_1,x_2) & \text{if $x_2=0$}, 
\\ 
z^{\gamma-3}(x_1+1,x_2-1) & \text{if $x_2\neq 0$}. 
\end{cases}
$}
\end{align}
Here $z^{\gamma}(x_1,x_2)$ denotes a one-soliton state labeled by the 
$U_q(A_1^{(1)})$-crystals. 
Since $\tilde{e}_2$ commutes with $T_r$, 
we can reduce the proof to the $U_q(A_1)$-highest weight elements 
\[
z^{\gamma_1}(l,0)\otimes z^{\gamma_2}(y_1,y_2)
\in \text{Aff}(\tilde{B}_l)\otimes \text{Aff}(\tilde{B}_k)
\quad \text{for $l>k$}. 
\]
From the second property of $T_{\natural}$, the proof is reduced to 
 the induction hypothesis on the number of applications of $T_{\natural}$.

\subsection{Plan of the paper}

This paper is organized as follows. 
In section 2 we recall the facts on the crystal and combinatorial $R$. 
We also recall a family of $U'_q(\mathfrak{g})$-crystals $\{B_l\}_{l\geq 1}$ 
for $\mathfrak{g}=A_n^{(1)}$ or $D_4^{(3)}$ and $U'_q(D_4^{(3)})$-crystal $B_{\natural}$ associated with the vertex $1$ and $2$ in the Dynkin diagram respectively. 
In section 3, we construct an algorithm for the combinatorial $R$ matrix 
in terms of Lecouvey's column insertion algorithm [17] for the tableaux of type $G_2$. 
In section 4, using the algorithm introduced in the previous section 
we construct a box ball system associated with $U'_q(D_4^{(3)})$-crystals 
and describe the scattering rule of solitons appearing in it.


\section*{Acknowledgement}
The author would like to thank Tetsuji Tokihiro for warm 
encouragements during the study. 
The author is also grateful to Atsuo Kuniba and Masato Okado 
for useful discussions. 
This research is partially supported by the 21 century COE program
at Graduate School of Mathematical Sciences, the University of Tokyo.

\section{Preliminaries}

\subsection{Crystals}\label{sec:crystal}

A crystal $B$ is a set $B$ with the map
$\tilde{e}_i, \tilde{f}_i : B\sqcup \{0\}\longrightarrow B\sqcup \{0\}$, 
satisfying the following properties: 
\begin{align}
& 
\text{$\tilde{e}_i 0=\tilde{f}_i 0=0$.}
\\ & 
\text{for any $b$ and $i$ there exist $n>0$ such that 
$\tilde{e}_i^nb=\tilde{f}_i^nb=0$.} 
\\ & 
\text{for $b, b'\in B$ and $i$, 
$\tilde{f}_i b=b' \iff b=\tilde{e}_i b'$.}
\end{align}
For $b\in B$, set 
\begin{align*}
\varepsilon_i(b)=\max\{n\geq 0 \ \vert \ \tilde{e}_i^n b\neq 0 \},\quad  
\varphi_i(b)=\max\{n\geq 0 \ \vert \ \tilde{f}_i^n b\neq 0 \}. 
\end{align*}
The crystal $B$ is identified with a colored oriented graph, 
called the crystal graph, if one draws an arrow as 
\begin{equation}
b\xrightarrow{i} b' \ \iff \ \tilde{f}_i b=b'. 
\end{equation}

For two crystals $B$ and $B'$, the tensor product $B\otimes B'$ of  
is defined by as a set 
$B\otimes B'=\left\{b\otimes b' \ \vert \ b\in B,\ b'\in B' \right\}$ 
with the actions of $\tilde{e}_i, \tilde{f}_i$ as follows: 
\begin{align*}
& 
{\tilde e}_i (b \otimes b')=
\begin{cases}
{\tilde e}_i b \otimes b' & \mbox{if}\; \varphi _i(b)\geq \varepsilon_i(b'), 
\\ 
b\otimes {\tilde e}_i b' & \mbox{if}\; \varphi _i(b)<\varepsilon_i(b').  
\end{cases}
\\ & 
{\tilde f}_i (b\otimes b')=
\begin{cases}
{\tilde f}_i b\otimes b' & \mbox{if}\;\varphi _i(b)>\varepsilon_i(b'), 
\\
b\otimes {\tilde f}_i b' & \mbox{if}\;\varphi_i(b) \leq \varepsilon_i(b'). 
\end{cases}
\end{align*}
Here $b\otimes 0$ and $0\otimes b'$ are understood to be 0.

Let $\mathfrak{g}$ be an affine Lie algebra 
and let $U'_q(\mathfrak{g})$ be the corresponding quantum affine algebra 
with out degree operator [14]. 
For $U'_q(\mathfrak{g})$-crystals $B$ and $B'$, assume that 
$B$ (resp. $B'$) is a finite set and $B\otimes B'$ is connected. 
Let $\text{Aff}(B)=\{z^d b \ \vert \ b\in B,\ d\in\mathbb{Z}\}$ be the affinization of a $U'_q(\mathfrak{g})$-crystal $B$. 
Then $\text{Aff}(B)$ is an infinite set, and is endowed with 
$\tilde{e}_i, \tilde{f}_i$ by 
$\tilde{e}_i(z^d b)=z^{d+\delta_{i0}}(\tilde{e}_i b)$,  
$\tilde{f}_i(z^d b)=z^{d-\delta_{i0}}(\tilde{f}_i b)$.

The combinatorial $R$-matrix (combinatorial $R$ for short) 
is a map with the following form: 
\begin{eqnarray*}
R:\ \mathrm{Aff}(B)\otimes \mathrm{Aff}(B') & \longrightarrow &
\mathrm{Aff}(B')\otimes \mathrm{Aff}(B)
\\
z^d b\otimes z^{d'} b' & \mapsto & 
z^{d'+H(b\otimes b')} \tilde{b'}\otimes z^{d-H(b\otimes b')}\tilde{b}, 
\end{eqnarray*}
{} Here a map $\iota: b\otimes b'\mapsto \tilde{b}'\otimes \tilde{b}$ is a crystal isomorphism, manely it commutes with $\tilde{e}_i, \tilde{f}_i$. 
$H(b\otimes b')$ is a $\mathbb{Z}$-valued function on $B\otimes B'$ 
satisfying the following property:

for any $b\otimes b'\in B\otimes B'$ and $i\in I$ such that 
$\tilde{e}_i(b\otimes b')=0$, 
\begin{equation*}
H({\tilde e}_i (b\otimes b'))=
\begin{cases}
H(b\otimes b')+1 &\text{if $i=0$, $\varphi_0(b)\geq \varepsilon_0(b')$, 
$\varphi_0(\tilde{b'})\geq \varepsilon_0(\tilde{b})$ },
\\
H(b\otimes b')-1 & \text{if $i=0$, $\varphi_0(b)<\varepsilon_0(b')$, 
$\varphi_0(\tilde{b'})<\varepsilon_0(\tilde{b})$},
\\
H(b\otimes b') & \text{otherwise}. 
\end{cases}
\end{equation*}
{} $H(b\otimes b')$ is called the energy function, and is 
determined up to a global additive constant. 
The existence of the crystal isomorphism 
$\iota: B\otimes B'\rightarrow B'\otimes B$ and energy function 
$H(b\otimes b')$ are guaranteed by the existence of the $R$ matrix in [14]. 
The definition of energy function assures that 
the combinatorial $R$ commutes with $\tilde{e}_i, \tilde{f}_i$.

Denote the crystal isomorphism 
\[
B\otimes B'\simeq B'\otimes B, \quad 
b\otimes b'\mapsto \tilde{b}'\otimes \tilde{b}
\] 
by the  following symbol. 
\[
\batten{b}{b'}{\tilde{b}'}{\tilde{b}}{a}
\]

\subsection{$U'_q(A_n^{(1)})$-crystal $B_l$}\label{sec:An crystal}

Let $\{B_l\}_{l\geq 1}$ be a family of crystals of the $l$-fold symmetric 
fusion of the vector representation of $U'_q(A_n^{(1)})$. 
It is as a set 
\begin{equation*}
B_l=\left\{ b=(x_1,x_2,\ldots,x_{n+1}) \ \vert \ x_i\in \mathbb{Z}_{\geq 0},\ \sum_{i=1}^n x_i=l \right\}. 
\end{equation*}
The actions $\tilde{e}_i, \tilde{f}_i$ ($i=0,1,\ldots ,n$) are given by 
\begin{align*}
\tilde{e}_0 b=&(x_1-1,\ldots,x_{n+1}+1),
\\ 
\tilde{e}_i b=&(x_1,\ldots ,x_i+1,x_{i+1}-1,\ldots,x_{n+1}) 
\quad \text{if $i\neq 0$}, 
\\ 
\tilde{f}_0 b=&(x_1+1,\ldots ,\ldots,x_{n+1}-1), 
\\
\tilde{f}_i b=&(x_1,\ldots ,x_i-1,x_{i+1}+1,\ldots,x_{n+1}) 
\quad \text{if $i\neq 0$}. 
\end{align*}

Let $\text{Aff}(B_l)$ be the affinization of $B_l$. 
A piecewise linear formula [3] for the combinatorial $R$ matrix 
for a crystal $B_l\otimes B_{l'}$
\begin{eqnarray*}
R: \text{Aff}(B_l)\otimes \text{Aff}(B_{l'}) & \longrightarrow & 
\text{Aff}(B_{l'})\otimes \text{Aff}(B_l)
\\ 
z^{d} b\otimes z^{d'} b' & \mapsto & 
z^{d'+H(b\otimes b')}\tilde{b}' \otimes z^{d-H(b\otimes b')}\tilde{b}
\end{eqnarray*}
is given as follows: 

for $b\otimes b'=(x_1,x_2,\ldots,x_{n+1})\otimes (y_1,y_2,\ldots,y_{n+1})
\in B_l \otimes B_{l'}$, 
\begin{align*}
& 
Q_i(b, b')=\min_{1\leq k\leq n+1}\left\{
\sum_{j=1}^{k-1}x_{i+j}+\sum_{j=k+1}^{n+1} y_{i+j} \right\}, 
\\ &
\tilde{b}'\otimes \tilde{b}=
(y'_1,y'_2,\ldots,y'_{n+1})\otimes (x'_1,x'_2,\ldots,x'_{n+1}), 
\\ & 
x'_i=x_i+Q_i(b, b')-Q_{i-1}(b, b'),\ y'_i=y_i+Q_{i-1}(b, b')-Q_i(b, b'), 
\\ &
H(b_1\otimes b_2)=-Q_{n+1}(b, b'), 
\end{align*}
where $x_{n+1+j}$ (resp. $y_{n+1+j}$) 
denotes $x_j$ (resp. $y_j$) for $1\leq j\leq n+1$.  

Another algorithm is known in [18], which is described by the tableau representation of $B_l$. Then, the element $b\in B_l$ is described by the Young tableau of shape $(l)$ with entries $1,2,\ldots , n+1$. 
\begin{equation}
(x_1,x_2,\ldots ,x_{n+1})\in B_l \iff 
(\underbrace{11...11}_{x_1}\underbrace{22...22}_{x_2}\ \ldots \ \underbrace{n+1...n+1}_{x_{n+1}})
\end{equation}
Here the frame of tableau is replaced by $(\ \cdots \ )$. 

\begin{example}
The $n=3$, $l=4$, $l'=2$ case. 

$\bullet$ Coordinate representation. 
\begin{align*}
& 
z^{d}(2,1,1,0)\otimes z^{d'}(0,1,1,0)\mapsto z^{d'+0}(1,1,0,0)\otimes z^{d-0}(1,1,2,0),  
\\ & 
z^{d}(2,1,1,0)\otimes z^{d'}(1,1,0,0)\mapsto z^{d'+1}(1,0,1,0)\otimes z^{d-1}(2,2,0,0),
\\ & 
z^{d}(0,1,1,2)\otimes z^{d'}(1,1,0,0)\mapsto z^{d'+2}(0,0,0,2)\otimes z^{d-2}(1,2,1,0).  
\end{align*}

$\bullet$ Tableau representation. 
\begin{align*} 
& 
z^{d}(1123)\otimes z^{d'}(23)\mapsto z^{d'+1}(12)\otimes z^{d-1}(1233), 
\\ & 
z^{d}(1123)\otimes z^{d'}(12)\mapsto z^{d'+1}(13)\otimes z^{d-1}(1122), 
\\ & 
z^{d}(2344)\otimes z^{d'}(12)\mapsto z^{d'+2}(44)\otimes z^{d-2}(1223). 
\end{align*}
\end{example}

Now we can state 
\begin{proposition}
(Yang-Baxter equation).\ 
The following relation holds on 
$\text{Aff}(B_l)\otimes \text{Aff}(B_{l'})\otimes \text{Aff}(B_{l''}):$ 
\begin{equation}
(R\otimes 1)(1\otimes R)(R\otimes 1)=(1\otimes R)(R\otimes 1)(1\otimes R). 
\end{equation}
\end{proposition}

\subsection{$U'_q(D_4^{(3)})$-crystal $B_l$}\label{sec:D_4^(3)}

A family of $U'_q(\mathfrak{g})$-crystals $\{B_l\}_{l\geq 1}$ 
corresponding to the vertex $1$ in the Dynkin diagram in 
the exceptional affine Lie algebra $\mathfrak{g}=D_4^{(3)}$ 
was studied in [16]. 
In this paragraph 
let us recall the coordinate representation of the crystal $B_l$.

As a set crystal $B_l$ is given by 
\begin{equation}
B_l=\bigoplus_{j=0}^l B(j\Lambda_1) \quad \text{as $U_q(G_2)$-crystals.}
\end{equation}

$B(j\Lambda_1)$ is isomorphic to the crystal base for finite-dimensional irreducible $U_q(G_2)$-module with highest weight $j\La_1$ in [12].  
Any element of $B(j\La_1)$ is represented as a one-row semistandard tableau 
with entries $1,2,3,0,\bar{3},\bar{2},\bar{1}$ as follows. 
\begin{equation}
\underbrace{
\begin{array}{c}
\hline
\multicolumn{1}{|c|}{\lw{1}\ldots \lw{1}} \\
\hline 
\end{array}
}_{w_1}\!
\underbrace{
\begin{array}{c}
\hline
\multicolumn{1}{|c|}{\lw{2}\ldots \lw{2}} \\
\hline 
\end{array}
}_{w_2}\!
\underbrace{
\begin{array}{c}
\hline
\multicolumn{1}{|c|}{\lw{3}\ldots \lw{3}} \\
\hline 
\end{array}
}_{w_3}\!
\underbrace{
\begin{array}{c}
\hline
\multicolumn{1}{|c|}{\ \lw{0}\ } \\
\hline 
\end{array}
}_{w_0}\!
\underbrace{
\begin{array}{c}
\hline
\multicolumn{1}{|c|}{\lw{\bar 3}\ldots \lw{\bar 3}} \\\hline 
\end{array}
}_{{\bar w}_3}\!
\underbrace{
\begin{array}{c}
\hline
\multicolumn{1}{|c|}{\lw{\bar 2}\ldots \lw{\bar 2}} \\
\hline 
\end{array}
}_{{\bar w}_2}\!
\underbrace{
\begin{array}{c}
\hline 
\multicolumn{1}{|c|}{\lw{\bar 1}\ldots \lw{\bar 1}} \\
\hline 
\end{array}
}_{{\bar w}_1}
\end{equation}
\begin{equation}
w_0\in \{0,1\}, \quad \sum_{i=1}^3(w_i+\bar{w}_i)+w_0=j.
\end{equation}
Now we set 
\begin{equation}
(x_1,x_2,x_3,{\bar x}_3,{\bar x}_2,{\bar x}_1)=
(w_1, w_2, 2w_3+w_0, 2\bar{w}_3+w_0, \bar{w}_2, \bar{w}_1). 
\end{equation}
Define a {\em coordinate representation} of $B(j\Lambda_1)$ by 
\[
B(j\Lambda_1)=\left\{
b=(x_1,x_2,x_3,{\bar x}_3,{\bar x}_2,{\bar x}_1)
\ \vert \ 
\begin{array}{l}
x_i, \bar{x}_i\in \mathbb{Z}_{\geq 0},
\\ 
x_3\equiv\bar{x}_3\;(\text{mod }2), \ 
\\ 
x_1+x_2+(x_3+{\bar x}_3)/2+{\bar x}_2+{\bar x}_1=j
\end{array}
\right\}. 
\]
Then the crystal structure is given as follows: 
\begin{align*}
{\tilde e}_0 b=&
\begin{cases}
(x_1 -1,\ldots ) &
\text{if $z_1+(z_2+(3z_4+(z_3)_-)_-)_-< 0$},
\\
(\ldots ,x_3 -1,\bar{x}_3 -1,\ldots ,\bar{x}_1 +1) &
\text{if $z_2+(3z_4+(z_1+z_3)_-)_-< 0\leq z_1$},
\\
(\ldots ,x_3 -2,\ldots ,\bar{x}_2 +1,\ldots ) &  
\text{if $3z_4+(z_3+(z_1)_-)_-< 0\leq z_2+(z_1)_+$},
\\
(\ldots ,x_2 -1,\ldots ,\bar{x}_3 +2,\ldots ) & 
\text{if $3z_4+(z_2+(z_1)_+)_+ > 0\geq z_3+(z_1)_+$},
\\
(x_1-1,\ldots ,x_3+1,\bar{x}_3+1,\ldots ) & 
\text{if $z_3+(3z_4+(z_1+z_2)_+)_+\geq 0> z_1$},
\\ 
(\ldots ,\bar{x}_1 +1) &  
\text{if $z_1+(z_3+(3z_4+(z_2)_+)_+)_+\geq 0$}.
\end{cases}
\\  
{\tilde f}_0 b=&
\begin{cases}
(x_1 +1,\ldots ) &  
\text{if $z_1+(z_2+(3z_4+(z_3)_-)_-)_-\leq 0$},
\\ 
(\ldots ,x_3 +1,\bar{x}_3 +1,\ldots ,\bar{x}_1 -1) &
\text{if $z_2+(3z_4+(z_1+z_3)_-)_-\leq 0<z_1$},
\\ 
(\ldots ,x_3 +2,\ldots ,\bar{x}_2 -1,\ldots ) &
\text{if $3z_4+(z_3+(z_1)_-)_-\leq 0<z_2+(z_1)_+$},
\\ 
(\ldots ,x_2 +1,\ldots ,\bar{x}_3 -2,\ldots ) & 
\text{if $3z_4+(z_2+(z_1)_+)_+ \geq 0>z_3+(z_1)_+$},
\\ 
(x_1 +1,\ldots ,x_3 -1,{\bar x}_3 -1,\ldots ) &  
\text{if $z_3+(3z_4+(z_1+z_2)_+)_+>0\geq z_1$},
\\ 
(\ldots,\bar{x}_1 -1) &  
\text{if $z_1+(z_3+(3z_4+(z_2)_+)_+)_+>0$}.
\end{cases}
\end{align*}
\begin{align*}
\veps_0(b)=&l-s(b)+(z_1+(z_2+3z_4+(z_3+(z_1)_+)_+)_+)_+-(2z_1+z_2+z_3+3z_4),\\ 
\vphi_0(b)=&l-s(b)+(z_1+(z_2+3z_4+(z_3+(z_1)_+)_+)_+)_+. 
\end{align*}
\begin{align*}
{\tilde e}_1 b=&
\begin{cases}
(\ldots,{\bar x}_2 +1,{\bar x}_1 -1) 
& \text{if ${\bar x}_2 -{\bar x}_3 \geq (x_2 -x_3)_+$}, 
\\  
(\ldots,x_3 +1,{\bar x}_3 -1,\ldots) 
& \text{if ${\bar x}_2 -{\bar x}_3 <0\leq x_3 -x_2$}, 
\\ 
(x_1 +1,x_2 -1,\ldots) 
& \text{if $({\bar x}_2 -{\bar x}_3)_+ <x_2 -x_3$}. 
\end{cases}
\\
{\tilde f}_1 b=&
\begin{cases}
(x_1 -1,x_2 +1,\ldots) 
& \text{if $({\bar x}_2 -{\bar x}_3)_+ \leq x_2 -x_3$}, 
\\
(\ldots,x_3 -1,{\bar x}_3 +1,\ldots) 
& \text{if ${\bar x}_2 -{\bar x}_3 \leq 0<x_3 -x_2$}, 
\\
(\ldots,{\bar x}_2 -1,{\bar x}_1 +1) 
& \text {if ${\bar x}_2 -{\bar x}_3 >(x_2 -x_3)_+$}. 
\end{cases}
\\
{\tilde e}_2 b=&
\begin{cases}
(\ldots,{\bar x}_3 +2,{\bar x}_2 -1,\ldots) 
& \text{if ${\bar x}_3 \geq x_3$}, 
\\
(\ldots,x_2 +1,x_3 -2,\ldots) 
& \text{if ${\bar x}_3 <x_3$}. 
\end{cases}
\\
{\tilde f}_2 b=&
\begin{cases}
(\ldots,x_2 -1,x_3 +2,\ldots) 
& \text{if ${\bar x}_3 \leq x_3$}, 
\\
(\ldots,{\bar x}_3 -2,{\bar x}_2 +1,\ldots) 
& \text{if ${\bar x}_3 >x_3$}. 
\end{cases}
\end{align*}
\begin{align*}
\veps_1(b)=&{\bar x}_1+({\bar x}_3-{\bar x}_2+(x_2-x_3)_+)_+,
\quad 
\veps_2(b)={\bar x}_2+\frac{1}{2}(x_3-{\bar x}_3)_+,
\\
\varphi_1(b)=&x_1+(x_3-x_2+({\bar x}_2-{\bar x}_3)_+)_+,
\quad 
\varphi_2(b)=x_2+\frac{1}{2}({\bar x}_3-x_3)_+.
\end{align*}
Here $(x)_+=\max(x,0)$, $(x)_-=\min(x,0)$, 
\begin{align*} \label{z1-4}
& 
z_1={\bar x}_1-x_1, \ 
z_2={\bar x}_2 -{\bar x}_3, \ 
z_3=x_3-x_2, \ 
z_4=({\bar x}_3-x_3)/2, 
\\ & 
s(b)=x_1+x_2+\frac{x_3+\bar{x}_3}{2}+\bar{x}_2+\bar{x}_1. 
\end{align*}

For instance, as an element of $B^{G_2}(7\La_1)$
\[
\begin{array}{ccccccc}
\hline
\multicolumn{1}{|c}{\lw{1}} &
\multicolumn{1}{|c}{\lw{2}} &
\multicolumn{1}{|c}{\lw{2}} &
\multicolumn{1}{|c}{\lw{3}} &
\multicolumn{1}{|c}{\lw{0}} &
\multicolumn{1}{|c}{\lw{\bar 1}} &
\multicolumn{1}{|c|}{\lw{\bar 1}} \\
\hline 
\end{array}
\iff (1,2,3,1,0,2). 
\]

Let us denote the elements of $B_1$ by 
\begin{align*}
\fbox{$1$}&=(1,0,0,0,0,0),& 
\fbox{$2$}&=(0,1,0,0,0,0),& 
\fbox{$3$}&=(0,0,2,0,0,0),\\
\fbox{$0$}&=(0,0,1,1,0,0),&
\fbox{${\bar 3}$}&=(0,0,0,2,0,0),&    
\fbox{${\bar 2}$}&=(0,0,0,0,1,0),\\
\fbox{${\bar 1}$}&=(0,0,0,0,0,1), &
\phi&=(0,0,0,0,0,0).
\end{align*}
then, the crystal graph of $B_1$ is given as follows: 
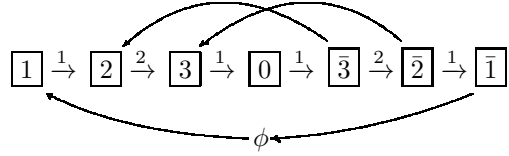
\begin{center}
\begin{picture}(220,80)
\put(0,30)
{\makebox(0,0){$\fbox{$1$}\xrightarrow{1}$}}
\put(30,30)
{\makebox(0,0){$\fbox{$2$}\xrightarrow{2}$}}
\put(60,30)
{\makebox(0,0){$\fbox{$3$}\xrightarrow{1}$}}
\put(90,30)
{\makebox(0,0){$\fbox{$0$}\xrightarrow{1}$}}
\put(120,30)
{\makebox(0,0){$\fbox{${\bar 3}$}\xrightarrow{2}$}}
\put(155,30)
{\makebox(0,0){$\fbox{${\bar 2}$}\xrightarrow{1}\fbox{${\bar 1}$}$}}
\put(82,3)
{\makebox(0,0){$\phi$}}
\bezier{200}(86,3)(125,5)(162,20)
\bezier{200}(1,20)(30,5)(77,3)
\bezier{200}(59,38)(100,70)(135,40)
\bezier{200}(30,38)(65,70)(107,40)
\bezier{200}(86,3)(87,3)(88,4)
\bezier{200}(86,3)(88,3)(89,2)
\bezier{200}(1,20)(2,19)(2,18)
\bezier{200}(1,20)(1,20)(3,20)
\bezier{200}(30,38)(31,39)(31,40)
\bezier{200}(30,38)(31,39)(32,39)
\bezier{200}(59,38)(60,39)(60,40)
\bezier{200}(59,38)(60,39)(62,39)
\end{picture}
\end{center}
The arrows without number are $0$-arrows.

\subsection{$U'_q(D_4^{(3)})$-crystal $B_{\natural}$}\label{sec:Bnatural}

It is known [6, 15, 16] that 
$U'_q(D_4^{(3)})$-module $V^2$ 
corresponding to the vertex $2$ in the Dynkin diagram 
has a crystal pseudobase and is decomposed into 
$V^2\simeq V(\overline{\Lambda}_2)\oplus V(\overline{\Lambda}_1)^{\oplus 2}\oplus V(0)$ 
as highest-weight $U_q(G_2)$-modules, where 
$\overline{\Lambda}_1=\Lambda_1-2\Lambda_0$ and 
$\overline{\Lambda}_2=\Lambda_2-3\Lambda_0$.  
The crystal $B_{\natural}$ of $V^2$ is given by as a set $B_{\natural}$ 
\begin{align*}
B_{\natural}=&
(\overline{\Lambda}_2)\oplus (\overline{\Lambda}_1)^{\oplus 2}\oplus B(0), 
\\  
B(0)=&\{\phi\}, 
\\ 
B(\overline{\Lambda}_1)^{\oplus 2}=&\bigsqcup_{j=1,2}\left\{
(\alpha)_j \ \vert \ \alpha=1,2,3,0,\bar{3},\bar{2},\bar{1} \right\}, 
\\ 
B(\overline{\Lambda}_2)=&
\left\{
\begin{pmatrix}
\alpha \\ \beta
\end{pmatrix}
\ \vert \ \text{The word $\alpha\beta$ is the element of $B(12)$.}
\right\}. 
\end{align*}
For a set of words $B(12)$, see section \ref{sec.lec}. 
The labeling of the elements of $B(\overline{\Lambda}_2)$ is taken from [12].Thus the actions $\tilde{e}_i, \tilde{f}_i$ ($i=1,2$) on 
$B(\overline{\Lambda}_2)$ agree with those in [12]. 
On $B(0)$ and $B(\overline{\Lambda}_1)^{\oplus 2}$ they are defined by 
\begin{align*}
& 
\tilde{e}_i\phi=\tilde{f}_i\phi=0 \quad (i=1,2), 
\\ & 
(1)_j\xrightarrow{1}
(2)_j\xrightarrow{2}
(3)_j\xrightarrow{1}
(0)_j\xrightarrow{1}
(\bar{3})_j\xrightarrow{2}
(\bar{2})_j\xrightarrow{1}
(\bar{1})_j 
\quad (j=1,2). 
\end{align*}
Here $b\xrightarrow{i}b'$ means $\tilde{f}_i b=b'$. 

The $0$-action on $B_{\natural}$ is given as follows. 
\begin{align*}
&
\begin{pmatrix}
\bar{3} \\ \bar{1}
\end{pmatrix}
\xrightarrow{0}(\bar{3})_1\xrightarrow{0}(2)_2\xrightarrow{0}
\begin{pmatrix}
1 \\ 2
\end{pmatrix}
,\quad 
\begin{pmatrix}
\bar{2} \\ \bar{1}
\end{pmatrix}
\xrightarrow{0}(\bar{2})_1\xrightarrow{0}(3)_2\xrightarrow{0}
\begin{pmatrix}
1 \\ 3
\end{pmatrix}
, \\ & 
\begin{pmatrix}
\bar{3} \\ \bar{2}
\end{pmatrix}
\xrightarrow{0}(0)_1\xrightarrow{0}(1)_2
,\quad 
(\bar{1})_1\xrightarrow{0}(0)_2\xrightarrow{0}
\begin{pmatrix}
2 \\ 3
\end{pmatrix}
,\quad 
(\bar{1})_2\xrightarrow{0}\phi \xrightarrow{0}(1)_1
,\\ & 
\begin{pmatrix}
0 \\ \bar{2}
\end{pmatrix}
\xrightarrow{0} (3)_1,
\quad 
(\bar{3})_2\xrightarrow{0}
\begin{pmatrix}
2 \\ 0
\end{pmatrix}
,\quad 
\begin{pmatrix}
0 \\ \bar{3}
\end{pmatrix}
\xrightarrow{0} (2)_1,
\quad 
(\bar{2})_2\xrightarrow{0}
\begin{pmatrix}
3 \\ 0
\end{pmatrix}
\\ & 
\begin{pmatrix} 0 \\ 0 \end{pmatrix},\quad 
\begin{pmatrix} 3 \\ \bar{3} \end{pmatrix},\quad 
\begin{pmatrix} 2 \\ \bar{3} \end{pmatrix},\quad 
\begin{pmatrix} 3 \\ \bar{2} \end{pmatrix}. 
\end{align*}

\subsection{Yang-Baxter equation}

For $U'_q(D_4^{(3)})$-crystals $B$ we define the affinization of them by 
\begin{equation}
\text{Aff}(B)=\left\{z^d b \ \vert \ d\in \mathbb{Z}, b\in B \right\}. 
\end{equation}
Here $z$ denotes spectral parameter. 

Define a combinatorial $R$ matrix for a $U_q(D_4^{(3)})$-crystals 
$B\otimes B'$ by 
\begin{eqnarray*}
R: \text{Aff}(B)\otimes \text{Aff}(B') & \longrightarrow & 
\text{Aff}(B')\otimes \text{Aff}(B)
\\ 
z^d b\otimes z^{d'} b' & \mapsto & 
z^{d'+H(b\otimes b')}\tilde{b} \otimes z^{d-H(b\otimes b')}\tilde{b}'
\end{eqnarray*}

Thanks to the ordinary Yang-Baxter equation, we obtain the following result.  

\begin{proposition}\label{prop:YBE}
(Yang-Baxter equation).\ 
The following equation holds on 
$\text{Aff}(B)\otimes \text{Aff}(B')\otimes \text{Aff}(B'')$: 
\begin{equation}
(R\otimes 1)(1\otimes R)(R\otimes 1)=(1\otimes R)(R\otimes 1)(1\otimes R). 
\end{equation}
for $B\otimes B'\otimes B''=B_l\otimes B_{l'}\otimes B_{l''}$ or 
$B_{\natural}\otimes B_l\otimes B_{l'}$. 
\end{proposition}

Hereafter we normalize the value of energy function on 
$B_l\otimes B_{l'}$ as follows. 
\begin{equation}\label{eq:nor.ene.}
H((l,0,0,0,0,0)\otimes (l',0,0,0,0,0))=0
\end{equation}

\section{Combinatorial $R$ matrix}\label{sec:comb.R}

\subsection{Lecouvey's column insertion algorithm for type $G_2$}\label{sec.lec}

We recall the column insertion for the $G_2$ case [17]. 
We define the set of words 
$B(1)$, $B(10)$, $B(12)$, $B(11)$, $B(112)$ and $B(121)$ as follows. 
Remark that $\alpha\beta\gamma$ denotes a word, and 
$\alpha \ \beta \ \gamma$ an one-row tableau. 

\begin{align*}
& 
B(1)=\{ 1,\ 2,\ 3,\ 0,\ \bar{3},\ \bar{2},\ \bar{1} \},
\ 
B(10)=
\{10,\ 1\bar{3},\ 1\bar{2},\ 2\bar{2},\ 2\bar{1},\ 3\bar{1},\ 0\bar{1} \},
\\ {} \\ &
B(12)=\{12,\ 13,\ 23,\ 20,\ 2\bar{3},\ 30,\ 3\bar{3},\ 3\bar{2},\ 00,\ 
0\bar{3},\ 0\bar{2},\ \bar{3}\bar{2},\ \bar{3}\bar{1},\ \bar{2}\bar{1}\},
\end{align*}
\[
B(11)=\left\{
\begin{array}{l} 
11,\ 21,\ 22,\ 31,\ 32,\ 33,\ 01,\ 02,\ 03,\ 
\bar{3}1,\ \bar{3}2,\ \bar{3}3,\ \bar{3}0,\ \bar{3}\bar{3},  
\\ 
\bar{2}1,\ \bar{2}2,\ \bar{2}3,\ \bar{2}0,\ \bar{2}\bar{3},\ \bar{2}\bar{2},\ 
\bar{1}1,\ \bar{1}2,\ \bar{1}3,\ \bar{1}0,\ \bar{1}\bar{3},\ \bar{1}\bar{2},\ 
\bar{1}\bar{1}
\end{array}
\right\},
\]

\[
B(112)=\left\{
\begin{array}{l}
112,\ 113,\ 212,\ 213,\ 223,\ 220,\ 22\bar{3},\ 312,\ 313,\ 323,\ 320,
\\ 
32\bar{3},\ 330,\ 33\bar{3},\ 33\bar{2},\ 012,\ 013,\ 023,\ 020,\ 
02\bar{3},\ 030,\ 03\bar{3},\\ 
03\bar{2},\ \bar{3}12,\ \bar{3}13,\ \bar{3}23,\ \bar{3}20,\ 
\bar{3}2\bar{3},\ \bar{3}30,\ \bar{3}3\bar{3},\ \bar{3}3\bar{2},\ \bar{3}00,\ 
\bar{3}0\bar{3}, 
\\ 
\bar{3}0\bar{2},\ \bar{3}\bar{3}\bar{2},\ \bar{3}\bar{3}\bar{1},\ 
\bar{2}12,\ \bar{2}13,\ \bar{2}23,\ \bar{2}20, 
\ \bar{2}2\bar{3},\ \bar{2}30,\ \bar{2}3\bar{3},\ \bar{2}3\bar{2},
\\ 
\bar{2}00,\ \bar{2}0\bar{3},\ 
\bar{2}0\bar{2},\ \bar{2}\bar{3}\bar{2},\ \bar{2}\bar{3}\bar{1},\ 
\bar{2}\bar{2}\bar{1}, \ \bar{1}12,\ \bar{1}13,\ \bar{1}23,\ 
\bar{1}20,\ \bar{1}2\bar{3},
\\ 
\bar{1}30,\ \bar{1}3\bar{3},\ \bar{1}3\bar{2},\ \bar{1}00,\ \bar{1}0\bar{3},\  
\bar{1}0\bar{2},\ \bar{1}\bar{3}\bar{2},\ \bar{1}\bar{3}\bar{1},\ 
\bar{1}\bar{2}\bar{1}
\end{array}
\right\},
\]

\[
B(121)=\left\{
\begin{array}{l}
121,\ 131,\ 122,\ 231,\ 201,\ 2\bar{3}1,\ 2\bar{3}2,\ 132,\ 133,\ 301,\ 
3\bar{3}1,
\\ 
3\bar{3}2,\ 3\bar{2}1,\ 3\bar{2}2,\ 3\bar{2}3,\ 
232,\ 233,\ 001,\ 0\bar{3}1,\ 0\bar{3}2,\ 0\bar{2}1,\ 0\bar{2}2,
\\ 
0\bar{2}3,\ 202,\ 203,\ 2\bar{3}3,\ 2\bar{3}0,\ 
2\bar{3}\bar{3},\ \bar{3}\bar{2}1,\ \bar{3}\bar{2}2,\ 
\bar{3}\bar{2}3,\ \bar{3}\bar{1}1,\ \bar{3}\bar{1}2,
\\ 
\bar{3}\bar{1}3,\ \bar{3}\bar{1}0,\ \bar{3}\bar{1}\bar{3},\ 
302,\ 303,\ 3\bar{3}3,\ 3\bar{3}0,\ 3\bar{3}\bar{3},\ 3\bar{2}0,\ 
3\bar{2}\bar{3},\ 3\bar{2}\bar{2},
\\ 
\bar{2}\bar{1}1,\ \bar{2}\bar{1}2,\ \bar{2}\bar{1}3,\ \bar{2}\bar{1}0,\ 
\bar{2}\bar{1}\bar{3},\ \bar{2}\bar{1}\bar{2},\ 002,\ 003,\ 0\bar{3}3,\ 
0\bar{3}0,\ 0\bar{3}\bar{3}, 
\\ 
0\bar{2}0,\ 0\bar{2}\bar{3},\ 0\bar{2}\bar{2},\ \bar{3}\bar{2}0,\ 
\bar{3}\bar{2}\bar{3},\ \bar{3}\bar{2}\bar{2},\ \bar{3}\bar{1}\bar{2},\ 
\bar{3}\bar{1}\bar{1},\ \bar{2}\bar{1}\bar{1} 
\end{array}
\right\}
\]

We also define the bijection 
$\xi: B(10)\rightarrow B(1)$ and $\eta: B(121)\rightarrow B(112)$ by 
from the $j$-th component in $B(10)$ (resp. $B(121)$) to that in 
$B(1)$ (resp. $B(112)$) respectively. 
e.g.,\ 
$\xi(2\bar{2})=0,\ \xi^{-1}(3)=1\bar{2}$, 
$\eta(0\bar{2}1)=030,\ \eta^{-1}(\bar{1}00)=\bar{3}\bar{2}0$.

Given a letter $\alpha\in \{1,2,3,0,\bar{3},\bar{2},\bar{1}, \phi\}$ 
and at most two row tableau $T$, which contains $\emptyset$, 
we define a tableau denoted by ($\alpha\rightarrow T$), 
and call such an algorithm a 
{\em column insertion of a letter $\alpha$ into a tableau $T$}. 
Let us begin with such $T$'s that have at most one column. 
The procedure of the column insertion $(\alpha \rightarrow T)$ 
is summarized as follows. 
\begin{align*}
& 
\text{case 0\quad 
$(\alpha \rightarrow \emptyset)=\alpha$}
\\ & 
\text{case 1\quad 
$
(\beta \rightarrow \alpha)=
\begin{array}{c}
\alpha \\ \beta 
\end{array}
$ 
\quad if $\alpha\beta \in B(12)$}
\\ & 
\text{case 2\quad 
$
(\beta\rightarrow \alpha)=(\beta\quad \alpha\rightarrow)$ 
\quad if $\alpha\beta \in B(11)$}
\\ &
\text{case 3 
$
\left(
\begin{array}{r}
\alpha \\ \gamma \rightarrow \beta
\end{array}
\right)=\left(
\begin{array}{l}
\alpha' \\ \beta' \quad \gamma'\rightarrow
\end{array}
\right); 
\ 
\gamma'\alpha'\beta'=\eta(\alpha \beta \gamma)
$
\quad if $\alpha \beta \gamma \in B(121)$}
\\ & 
\text{case 4\quad 
$(\beta\rightarrow \alpha)=(\ \gamma \rightarrow \ )$; 
$\gamma=\xi(\alpha\beta)$ \quad if $\alpha\beta \in B(10)$}
\\ & 
\text{case 5\quad 
$
(\beta\rightarrow \alpha)=
\begin{cases}
\beta & \text{if $\alpha=\phi$}, 
\\ 
\alpha & \text{if $\beta=\phi$}, 
\\ 
\emptyset & \text{if $(\alpha,\beta)=(\phi,\phi)$ or $(1,\bar{1})$}. 
\end{cases}
$}
\end{align*}
{} The above cases do not cover all the tableaux with two rows, 
but we only deal with these situations in this paper. 
When $T$ is a tableau with shape $(m,n)$ for $n\in \{0,1\}$, 
we repeat the above procedure: 
we insert a box into the leftmost column of $T$ according to the above formula.
If it is not case 0 or case 1, replace the column by the right hand side of the formula. Otherwise, replace the column by the right hand side of the formula without the right box. 
We say that this right box is bumped. 
We insert it into the second column of $T$ from the left and 
repeat the procedure above until we come to the case 0 or case 1.


\begin{example}\label{ex:bumping}
For any $n\in \mathbb{Z}_{>0}$, let 
$b_1\otimes b_2=(\alpha_1\ \alpha_2\ \ldots \ \alpha_n)\otimes (\beta)$ 
be the tableau representation of the elements in $B_l\otimes B_1$. 
Let $b_1*b_2$ be the tableau obtained by column insertion 
$(b_2\rightarrow b_1)$. Then we have 

\begin{align}
& 
\text{$b_1=(\phi)$ case:\quad 
$b_1*b_2=\emptyset$ if $\beta=\phi$; $=\beta$ if $\beta\neq \phi$.} 
\\ & 
\text{$b_1\neq (\phi)$, $b_2=(\phi)$ case:\quad 
$b_1*b_2=\alpha_1\ \alpha_2\ \ldots \ \alpha_n$.} 
\\ & 
\text{$\alpha_1\beta\in B(11)$ case:\quad 
$b_1*b_2=\beta\ \alpha_1\ \alpha_2\ \ldots \ \alpha_n$.} 
\\ & 
\text{$(\alpha_1,\beta)=(1,\bar{1})$ case:\quad 
$b_1*b_2=\alpha_2\ \alpha_3\ \ldots \ \alpha_n$.} 
\\ & 
\text{$\alpha_1\beta\in B(12)$ case:\quad 
$
b_1*b_2=
\begin{array}{l}
\alpha_1 \ \alpha_2 \ \alpha_3 \ \ldots \ \alpha_n \\ \beta
\end{array}
$.} 
\\ & 
\text{$\alpha_1\beta\in B(10)$, $\alpha_2\xi(\alpha_1\beta)\in B(11)$ 
case:\quad $b_1*b_2=\gamma \ \alpha_2\ \ldots \ \alpha_n $.} 
\\ & 
\text{$\alpha_1\beta\in B(10)$, $\alpha_2\xi(\alpha_1\beta)\in B(12)$ 
case:\quad 
$
b_1*b_2=
\begin{array}{l}
\alpha_2 \ \alpha_3 \ \alpha_4 \ \ldots \ \alpha_n \\ \gamma
\end{array}
$.} 
\end{align}
where $\gamma=\xi(\alpha_1\beta)$ if $\alpha_1\beta\in B(10)$. 
\end{example}

We define the reverse bumping algorithm by 
$\text{RHS}=\text{LHS}$ for the cases $0$-$3$.

\begin{example}\label{prop:reverse bumping}
Fix $n\in \mathbb{Z}_{\geq 1}$. 
The procedure of bumping out all letters from the tableau 
\[
T_n=(r'_n\rightarrow (r'_{n-1}\rightarrow (\cdots (r'_1\rightarrow 
(q_0 \rightarrow (p_0\rightarrow \emptyset)))\cdots )))
\]
is given by 
\begin{align*}
T_n=&
\begin{array}{l}
p_n \ r_n \ r_{n-1}\cdots r_2 \ (r_1\rightarrow \emptyset) \\ q_n
\end{array}
\\=&
\begin{array}{l}
p_n \ r_n \ r_{n-1}\cdots (r_{i+1}\rightarrow r_i) \cdots r_1 \\ q_n
\end{array}
\quad \text{for $1\leq i\leq n-1$}
\\=& 
\left(
\begin{array}{l}
p_n \\ q_n \quad \gamma_n \rightarrow
\end{array}
\right)r_{n-1}\cdots r_1
=
\left(
\begin{array}{r}
p_{n-1} \\ \gamma'_n \rightarrow q_{n-1}
\end{array}
\right)r_{n-1}\cdots r_1
\\=& 
(r'_n \rightarrow T_{n-1})=
(r'_n \rightarrow (r'_{n-1} \rightarrow T_{n-2}))=\cdots 
\\=& 
(r'_n \rightarrow (r'_{n-1} \rightarrow ( \cdots 
(r'_1 \rightarrow T_0) \cdots ))) \quad 
\\=& 
(r'_n\rightarrow (r'_{n-1}\rightarrow (\cdots (r'_1\rightarrow 
(q_0 \rightarrow (p_0\rightarrow \emptyset)))\cdots )))
\end{align*}
where \ $p_{i-1}q_{i-1}r'_i=\eta^{-1}(r_i p_i q_i)$ for $n\geq i \geq 1$.  
\end{example}

\subsection{Algorithm}\label{sec:algorithm for R matrix}

Let $\{B_l\}_{\geq 1}$ be the family of $U'_q(D_4^{(3)})$-crystals. 
In this paragraph we will give the explicit description of the combinatorial 
$R$ on $B_l\otimes B_1$. 
Fix the positive integer $l$. For a given element $b_1\otimes b_2$ in $B_l\otimes B_1$, we apply the following procedure.

\begin{algo}
Step 1.\  
Denote $b_1\otimes b_2\in B_l\otimes B_1$ by 
tableau representation 
\begin{align}
b_1\otimes b_2=(\alpha_1\ \alpha_2\ \ldots \ \alpha_n) \otimes (\beta)
\ \in B_l\otimes B_1 \quad (1\leq n\leq l). 
\end{align}

Step 2.\ 
Compute the column insertion $(b_2\rightarrow b_1)$, 
and obtain the tableau $b_1*b_2$. 
Let $(p+q,q)$ be the shape of $b_1*b_2$.   

Step 3.\ 
Define a map $H: B_l\otimes B_1\rightarrow \mathbb{Z}$ by 
\begin{align}
H(b_1\otimes b_2)=
\begin{cases}
-2 \quad \text{if 
$\alpha_1\beta\in B(10)$, $\alpha_2\gamma \in B(11)$}, 
\\ 
\max\left(-2,\ (b_1*b_2)_1-l-1\right)\quad \text{otherwise}. 
\end{cases}
\end{align}
where $(b_1*b_2)_1=p+q$. 

Step 4.\ 
Bump out all letters from $b_1*b_2$ by the reverse bumping algorithm. 
Denote the procedure as follows: 
\begin{align}
b_1*b_2=(t_1\rightarrow (t_2\rightarrow (\ldots 
(t_{p+2q}\rightarrow \emptyset)\ldots ))) ;\ t_0=\emptyset. 
\end{align}
Set $(\mathscr{T}_r)=(t_1\ t_2\ \ldots t_r);\ 0\leq r\leq p+2q$. 

Step 5.\ 
Define $b'_2\otimes b'_1\in B_1\otimes B_l$ by 
\begin{align}
b'_2\otimes b'_1=
\begin{cases}
(t_{p+2q})\otimes (\mathscr{T}_{p+2q-1}) & \text{for case I}, 
\\ 
(t'_{p+2q})\otimes (\mathscr{T}_{p+2q-1}t''_{p+2q}) & \text{for case II}, 
\\ 
(1)\otimes (\mathscr{T}_{p+2q}\bar{1}) & \text{for case III}, 
\\ 
(\mathscr{T}_{p+2q})\otimes (\phi) & \text{for case IV}, 
\\ 
(\phi)\otimes (\mathscr{T}_{p+2q}) & \text{for case V}, 
\\ 
(\phi)\otimes (\phi) & \text{otherwise}. 
\end{cases}
\end{align}
where 
\[
\begin{array}{rl}
III.& \alpha_1\beta\in B(11),\ n<l-1, 
\\ 
V.& \alpha_1\beta\in B(11),\ n=l-1, 
\\ 
I.& \alpha_1\beta\in B(11),\ n=l, 
\\ 
II.& \alpha_1\beta\in B(12),\ n<l;\ t'_{p+2q}t"_{p+2q}=\xi^{-1}(t_{p+2q}), 
\\ 
I.& \alpha_1\beta \in B(12),\ n=l, 
\\ 
I.& \alpha_1\beta\in B(10), 
\\ 
III.& (\alpha_1,\beta)=(1,\bar{1}),\ n=1, 
\\ 
IV.& (\alpha_1,\beta)=(1,\bar{1}),\ n=2, 
\\ 
I.& (\alpha_1,\beta)=(1,\bar{1}),\ n>2, 
\\ 
IV.& b_1\neq (\phi),\ b_2=(\phi),\ l=1, 
\\ 
V.& b_1\neq (\phi),\ b_2=(\phi),\ l\neq 1,\ n<l, 
\\ 
I.& b_1\neq (\phi),\ b_2=(\phi),\ l\neq 1,\ n=l,
\\ 
V.& b_1=(\phi),\ b_2\neq (\phi),\ l=1, 
\\ 
III.& b_1=(\phi),\ b_2\neq (\phi),\ l\neq 1. 
\end{array}
\]

Step 6.\ 
Define the image of a map $R: \text{Aff}(B_l)\otimes \text{Aff}(B_1)\rightarrow \text{Aff}(B_1)\otimes \text{Aff}(B_l)$ by 
\begin{equation}\label{eq:algo R}
R : z^{\gamma_1}b_1\otimes z^{\gamma_2}b_2
\mapsto z^{\gamma_2+H(b_1\otimes b_2)}b'_2\otimes z^{\gamma_1-H(b_1\otimes b_2)}b'_1. 
\end{equation}
\end{algo}

Now the theorem is 
\begin{theorem}\label{thm:Combinatorial R}
Assume (\ref{eq:nor.ene.}). 
Then the map (\ref{eq:algo R}) is the combinatorial $R$ matrix 
for $U'_q(D_4^{(3)})$-crystals. 
\end{theorem}

\begin{example}
The $l=3$ case.  
\[
R: z^{\gamma_1}(2\bar{2}\bar{1})\otimes z^{\gamma_2}(\bar{2})\mapsto 
z^{\gamma_2+(-2)}(\bar{1})\otimes z^{\gamma_1-(-2)}(0\bar{2}). 
\]
Because of the bumping algorithm and reverse bumping algorithm, we have 
\[
(2\bar{2}\bar{1})*(\bar{2})=(\bar{2}\rightarrow 2\bar{2}\bar{1})
=0\bar{2}\bar{1},\quad 
H((2\bar{2}\bar{1})\otimes (\bar{2}))=-2, 
\]
\[
0\bar{2}\bar{1}=(0\rightarrow \bar{2}\bar{1})=
(0\rightarrow (\bar{2}\rightarrow \bar{1})). 
\]
\end{example}

\begin{example}
The $l=4$ case.  
\[
R:\ z^{\gamma_1}(30\bar{1}) \otimes z^{\gamma_2}(\bar{3}) \mapsto 
z^{\gamma_2+(-2)}(2) \otimes z^{\gamma_1-(-2)}(2\bar{2}\bar{2}\bar{1}) 
\]
Because of the bumping algorithm and reverse bumping algorithm, we have 
\[
(30\bar{1})*(\bar{3})=
\begin{array}{l}
3 \ 0 \ \bar{1} \\ \bar{3}
\end{array}
,\ 
H((30\bar{1})\otimes (\bar{3}))=
\max\left(-2,((30\bar{1})*(\bar{3}))_1 -4-1\right)=-2, 
\]
\[
\begin{array}{l}
3 \ 0 \ \bar{1} \\ \bar{3}
\end{array}
=\left(
2\rightarrow 
\begin{array}{l}
0 \ \bar{1}  \\ \bar{2}
\end{array}
\right)=\left(
2\rightarrow \left(\bar{2}\rightarrow 
\begin{array}{l}
\bar{3}  \\ \bar{2}
\end{array}
\right) \right)=
(2\rightarrow (\bar{2}\rightarrow (\bar{2}\rightarrow (\bar{3}\rightarrow 
\emptyset )))). 
\]
Since $(\bar{3}\rightarrow \emptyset)$ is replaced to $(\bar{1}\rightarrow 2)$ 
in $l=4$ case, we have 
\[
\begin{array}{l}
3 \ 0 \ \bar{1} \\ \bar{3}
\end{array}
=(2\rightarrow (\bar{2}\rightarrow (\bar{2}\rightarrow 
(\bar{1}\rightarrow 2)))). 
\]
\end{example}

\subsection{Combinatorial $R$ for the highest weight elements}

We call an element $b$ of a $U'_q(D_4^{(3)})$ crystal a $U_q(G_2)$ 
highest weight element if it satisfies 
\begin{equation}
\tilde{e}_i b=0 \quad  \text{for $i=1,2$}.
\end{equation} 
Let $\hat{b}_1\otimes \hat{b}_2\in B_l\otimes B_1$ be a $U_q(G_2)$ 
highest weight element. Then, 
\begin{equation}
\hat{b}_1\otimes \hat{b}_2=
(n,0,0,0,0,0)\otimes (x_1,x_2,x_3,\bar{x}_3,0,\bar{x}_1)
\end{equation}
\begin{align}
(x_2,x_3,\bar{x}_3,\bar{x}_1)=(0,0,0,0),\ 
s(\hat{b}_2)\leq 1 & \qquad \text{for $n=0$}, 
\\  
(x_3, \bar{x}_3)=(1,1),\ (0,0),\ s(\hat{b}_2)\leq 1
 & \qquad \text{if $n=1$}, 
\\  
x_3\neq 2, \quad s(\hat{b}_2)\leq 1 
& \qquad \text{otherwise}. 
\end{align}

For any $m\in \mathbb{Z}_{\geq 0}$, 
we denote $(m,0,0,0,0,0)\in B_l$ by $(1^m)$. For example, 
\begin{align*}
& 
(0,0,0,0,0,0)=(\phi)=(1^0),
\\ &
(m,1,1,1,0,0)=(20) \quad \text{if $m=0$}
\ ; \ =(1^m20) \quad \text{otherwise}. 
\end{align*}

\begin{proposition}\label{lem:comb.R for h.w.e.}
Let a map $\iota: B_l\otimes B_1 \rightarrow B_1\otimes B_l$ be 
the crystal isomorphism for $U'_q(D_4^{(3)})$ crystals. 
For $0\leq n\leq l$ set 
$\hat{b}_1\otimes \hat{b}_2=(1^n)\otimes (\beta)\in B_l\otimes B_1$.   
Then we have 
\begin{equation}\label{eq:Rhwe}
\begin{array}{l}
\iota(\hat{b}_1\otimes \hat{b}_2)= 
\begin{cases}
(1)\otimes (1^{n+1}\bar{1}) & \text{if $\beta=1$, $0\leq n\leq l-2$},
\\ 
(\phi) \otimes (1^l) & \text{if $\beta=1$, $n=l-1$},
\\ 
(1) \otimes (1^l) & \text{if $\beta=1$, $n=l$}, 
\\
(1) \otimes (1^{n-1}20) & \text{if $\beta=2$, $1\leq n\leq l-1$},
\\ 
(1) \otimes (1^{l-1}2) & \text{if $\beta=2$, $n=l$}, 
\\ 
(1)\otimes (1^{n-1}0)
& \text{if $\beta=0$, $1\leq n\leq l$}, 
\\ 
(1) \otimes (1^{n-2}2) & \text{if $\beta=\bar{3}$, $2\leq n\leq l$}. 
\\ 
(1) \otimes (\bar{1}) & \text{if $\beta=\bar{1}$, $n=1$},
\\ 
(1) \otimes (\phi) & \text{if $\beta=\bar{1}$, $n=2$}, 
\\ 
(1) \otimes (1^{n-2}) & \text{if $\beta=\bar{1}$, $3\leq n\leq l$}. 
\\ 
(\phi) \otimes (1^n) & \text{if $\beta=\phi$, $0\leq n\leq l-1$}, 
\\ 
(1) \otimes (1^{l-1}) & \text{if $\beta=\phi$, $n=l$}. 
\end{cases}
\end{array}
\end{equation}
\end{proposition}
{} We will give a proof of this proposition in Appendix A. 

\begin{corollary}\label{cor:comb.R for h.w.e.}
The energy function associated with the crystal isomorphism 
$(\ref{eq:Rhwe})$ is given by 
\begin{equation}
H(\hat{b}_1\otimes \hat{b}_2)=
\begin{cases}
-2 \quad \text{if $\hat{b}_1\otimes \hat{b}_2=(1^n)\otimes (0)$}, 
\\ 
\max(-2, (\hat{b}_1*\hat{b}_2)_1-l-1) \quad \text{otherwise}. 
\end{cases}
\end{equation}
\end{corollary}
\begin{proof}
Thanks to Example \ref{ex:bumping}, we have 
\[
(\hat{b}_1*\hat{b}_2)_1=((1^n)*(\beta))_1=
\begin{cases}
0 & \text{if $n=0$, $\beta=\phi$}, 
\\ 
1 & \text{if $n=0$, $\beta\neq \phi$}, 
\\ 
n+1 & \text{if $n\neq 0$, $\beta=1$}, 
\\ 
n & \text{if $n\neq 0$, $\beta=2,0,\phi$}, 
\\ 
n-1 & \text{if $n\neq 0$, $\beta=\bar{3}, \bar{1}$}. 
\end{cases}
\]
Thanks to Appendix A, 
the energy function $H(\hat{b}_1\otimes \hat{b}_2)$ is given as follows. 
\begin{align*}
H(\hat{b}_1\otimes \hat{b}_2)=&
\begin{cases}
0 & \text{if $l=1$, $\hat{b}_1\otimes \hat{b}_2=(1)\otimes (1)$},
\\ 
-1 & 
\text{if $l=1$, 
$\hat{b}_1\otimes \hat{b}_2=(1)\otimes (2),\ 
(1)\otimes (\phi),\ (\phi)\otimes (1)$},
\\
0 & \text{if $l=2$, $\hat{b}_1\otimes \hat{b}_2=(11)\otimes (1)$},
\\ 
-1 & \text{if 
$l=2$, $\hat{b}_1\otimes \hat{b}_2=
(11)\otimes (2),\ (11)\otimes (\phi),\ (1)\otimes (1)$}, 
\\
0 & \text{if $l\geq 3$, $\hat{b}_1\otimes \hat{b}_2=(1^l)\otimes (1)$}, 
\\ 
-1 & \text{if $l\geq 3$, 
$\hat{b}_1\otimes \hat{b}_2=
(1^l)\otimes (2),\ (1^l)\otimes (\phi),\ (1^{l-1})\otimes (1)$}, 
\\ 
-2 & \text{otherwise}. 
\end{cases}
\end{align*}

Hence the statement is valid. 
\end{proof}

\subsection{Column insertion and $U_q(G_2)$ crystal morphism}

For a dominant integral weight $\lambda$, 
let $B(\lambda)$ be the $U_q(G_2)$ crystal associated with 
the irreducible highest weight representation $V(\lambda)$ [12]. 
The elements of $B(\lambda)$ can be represented by semistandard tableaux for 
$G_2$ of shape $\lambda$.  

Let $B(\lambda)\otimes B(\lambda')\simeq \bigoplus_j B(\lambda_j)^{\oplus m_j}$ be the tensor product decomposition of crystals, 
where $\lambda_j$'s are distinct highest weights and 
$m_j (\geq 1)$ is multiplicity of $B(\lambda_j)$. 
Forgetting the multiplicities we have the canonical morphism from 
$B(\lambda)\otimes B(\lambda')$ to $\bigoplus_j B(\lambda_j)$. 
The following proposition is due to [17]. 

\begin{proposition}\label{prop:G2 morphism}
Let $b_1*b_2$ be the tableau obtained from successive column insertions of 
letters of the Japanese reading word of $b_2$ into $b_1$. 
Define a map 
$\Phi: B(\lambda)\otimes B(\lambda') \longrightarrow \bigoplus_j B(\lambda_j)$
 by $\Phi(b_1\otimes b_2)=b_1*b_2$.
Then, the map $\Phi$ gives the unique crystal morphism. 
\end{proposition}

Thanks to this proposition, we have 
\begin{lemma}\label{lem:G2}
Let $\iota: B_l\otimes B_1\rightarrow B_1\otimes B_l$ 
be the crystal isomorphism for $U'_q(D_4^{(3)})$-crystal. 
Let $\hat{b}_1\otimes \hat{b}_2\in B_l\otimes B_1$ and 
$\hat{b}''_2 \otimes \hat{b}''_1\in B_1\otimes B_l$ 
be the $U_q(G_2)$-highest weight elements. 
For given $b_1\otimes b_2\in B_l\otimes B_1$ and 
$b''_2\otimes b''_1\in B_1\otimes B_l$, assume that  
\begin{itemize}
\item[$(1)$] \text{$b_1*b_2=b''_2*b''_1$.}

\item[$(2)$] \text{$b_1\otimes b_2$ (resp. $b''_2*b''_1$) is connected to 
$\hat{b}_1\otimes \hat{b}_2$ (resp. $\hat{b}''_2 \otimes \hat{b}''_1$) 
by $1$- and $2$-arrows.}

\item[$(3)$] \text{$\iota(\hat{b}_1\otimes \hat{b}_2)=\hat{b}''_2\otimes \hat{b}''_1$.} 
\end{itemize}
Then we have $\iota(b_1\otimes b_2)=b''_2\otimes b''_1$. 
\end{lemma}

\subsection{Proof of Theorem \ref{thm:Combinatorial R}}

First, we show the statement for the energy function. 
Since $\tilde{f}_1$ and $\tilde{f}_2$ does not change the shape of 
the tableaux [12] and $H(b_1\otimes b_2)=H(\hat{b}_1\otimes \hat{b}_2)$, 
the algorithm for energy function is reduced to 
Corollary \ref{cor:comb.R for h.w.e.}. 
Hence we have 
\[
H(b_1\otimes b_2)=H(\hat{b}_1\otimes \hat{b}_2)=
\begin{cases}
-2 \quad \text{if $\alpha_1\beta\in B(10)$, $\alpha_2\gamma\in B(11)$}, 
\\ 
\max(-2, (b_1*b_2)_1-l-1) \quad \text{otherwise}. 
\end{cases}
\]

Next, we will show the statement for the crystal isomorphism. 
By definition, the explicit form of $b'_2\otimes b'_1\in B_1\otimes B_l$ 
is given by 
\[
b'_2\otimes b'_1=
\begin{cases}
(1) \otimes (\beta\ \alpha_1\ \ldots \ \alpha_n \ \bar{1}) & 
\text{if $\alpha_1\beta\in B(11)$, $n< l-1$}, 
\\ 
(\phi) \otimes (\beta \ \alpha_1 \ \ldots \ \alpha_{l-1}) & 
\text{if $\alpha_1\beta\in B(11)$, $n=l-1$}, 
\\ 
(\alpha_l) \otimes (\beta \ \alpha_1 \ \ldots \ \alpha_{l-1}) & 
\text{if $\alpha_1\beta\in B(11)$, $n=l$}, 
\\ 
(1)\otimes (\bar{1}) & 
\text{if $(\alpha_1,\beta)=(1,\bar{1})$, $n=1$}, 
\\ 
(\alpha_2) \otimes (\phi) & 
\text{if $(\alpha_1,\beta)=(1,\bar{1})$, $n=2$}, 
\\ 
(\alpha_n) \otimes (\alpha_2 \ \alpha_3 \ \ldots \ \alpha_{n-1}) 
& \text{if $(\alpha_1,\beta)=(1,\bar{1})$, $n>2$}, 
\\ 
(\alpha_n) \otimes (\gamma \ \alpha_2 \ \alpha_3 \ \ldots \ \alpha_{n-1}) 
& \text{if 
$\alpha_1\beta\in B(10)$, $\alpha_2\gamma\in B(11)$}, 
\\ 
(p_0) \otimes (r'_{n-2} \ r'_{n-3} \ \ldots \ r'_1 \ q_0) 
& \text{if 
$\alpha_1\beta\in B(10)$, $\alpha_2\gamma\in B(12)$}, 
\\ 
(p'_0) \otimes (r'_{n-1} \ r'_{n-2} \ \ldots \ r'_1 \ q_0 \ p''_0) & 
\text{if $\alpha_1\beta\in B(12)$, $n<l$}, \\ 
(p_0) \otimes (r'_{l-1}\ r'_{l-2}\ \ldots \ r'_1\ q_0) & 
\text{if $\alpha_1\beta\in B(12)$, $n=l$}, 
\\ 
(\alpha_1) \otimes (\phi) & 
\text{if $b_1 \neq (\phi)$, $b_2=(\phi)$, $l=1$}, 
\\ 
(\phi) \otimes (\alpha_1 \ \alpha_2 \ \ldots \ \alpha_n) & 
\text{if $b_1\neq (\phi)$, $b_2=(\phi)$, $l\neq 1$, $n<l$}, 
\\ 
(\alpha_l) \otimes (\alpha_1 \ \alpha_2 \ \ldots \ \alpha_{l-1}) & 
\text{if $b_1\neq (\phi)$, $b_2=(\phi)$, $l\neq 1$, $n=l$}, 
\\ 
(\phi) \otimes (\beta) & \text{if $b_1=(\phi)$, $b_2\neq (\phi)$, $l=1$}, 
\\ 
(1) \otimes (\beta \ \bar{1}) 
& \text{if $b_1=(\phi)$, $b_2\neq (\phi)$, $l\neq 1$}, 
\\ 
(\phi) \otimes (\phi) & \text{if $b_1=(\phi)$, $b_2=(\phi)$}. 
\end{cases}
\]
where the new letters are defined as follows. 
\begin{itemize}
\item[$(i)$] 
The $\alpha_1\beta\in B(10)$ case. \quad Put 
\[
\gamma=\xi(\alpha_1\beta). 
\]
\item[$(ii)$] 
The $\alpha_1\beta\in B(10)$, $\alpha_2\gamma\in B(12)$ case.\quad Put 
\[
\left\{
\begin{array}{l}
(r_1,r_2,\ldots ,r_{n-2},p_{n-2},q_{n-2})=
(\alpha_n,\alpha_{n-1},\ldots ,\alpha_3,\alpha_2,\gamma), 
\\  
p_{i-1}q_{i-1}r'_i=\eta^{-1}(r_ip_iq_i) \quad \text{for $n-2\geq i\geq 1$}. 
\end{array}
\right.
\]
\item[$(iii)$] 
The $\alpha_1\beta\in B(12)$, $n=l$ case. \quad  Put 
\[
\left\{
\begin{array}{l}
(r_1,r_2,\ldots ,r_{l-1},p_{l-1},q_{l-1})=
(\alpha_l,\alpha_{l-1},\ldots ,\alpha_2,\alpha_1,\beta), 
\\ 
p_{i-1}q_{i-1}r'_i=\eta^{-1}(r_ip_iq_i) \quad \text{for $l-1\geq i\geq 1$}. 
\end{array}
\right.
\]
\item[$(iv)$] The $\alpha_1\beta\in B(12)$, $n<l$ case. \quad  Put 
\[
\left\{
\begin{array}{l}
(r_1,r_2,\ldots ,r_{n-1},p_{n-1},q_{n-1})=
(\alpha_n,\alpha_{n-1},\ldots ,\alpha_2,\alpha_1,\beta), 
\\ 
p_{i-1}q_{i-1}r'_i=\eta^{-1}(r_ip_iq_i) \quad \text{for $n-1\geq i\geq 1$}, 
\\ 
p'_0p''_0=\xi^{-1}(p_0). 
\end{array}
\right.
\]
\end{itemize}

Let $\hat{b}_1\otimes \hat{b}_2$ (resp. $\hat{b}'_2\otimes \hat{b}'_1$) 
be the $U_q(G_2)$ highest weight elements  
which are connected to $b_1\otimes b_2$ (resp. $b'_2\otimes b'_1$) 
by $1$- and $2$-arrows. 
By direct computations, we have the following equation. 

\begin{align*}
\hat{b}_1\otimes \hat{b}_2=&
\begin{cases}
(1^n)\otimes (1) & \text{if $\alpha_1\beta\in B(11)$, $1\leq n\leq l$}, 
\\ 
(1^n)\otimes (\bar{1}) 
& \text{if $(\alpha_1,\beta)=(1,\bar{1})$, $1\leq n\leq l$},
\\ 
(1^n) \otimes (0) 
& \text{if $\alpha_1\beta\in B(10)$, $\alpha_2\gamma\in B(11)$, 
$1\leq n\leq l$},
\\ 
(1^n) \otimes (\bar{3}) 
& \text{if $\alpha_1\beta\in B(10)$, $\alpha_2\gamma\in B(12)$, 
$2\leq n\leq l$}, 
\\ 
(1^n) \otimes (2) 
& \text{if $\alpha_1\beta\in B(12)$}, 
\\ 
(1^n) \otimes (\phi) & \text{if $b_1\neq (\phi)$, $b_2=(\phi)$, 
$1\leq n\leq l$}, 
\\ 
(\phi)\otimes (1) & \text{if $b_1=(\phi)$, $b_2\neq (\phi)$}, 
\\ 
(\phi)\otimes (\phi) & \text{if $b_1=(\phi)$, $b_2=(\phi)$}. 
\end{cases}
\\ 
\hat{b}'_2\otimes \hat{b}'_1=&
\begin{cases}
(1)\otimes (1^{n+1}\bar{1}) 
& \text{if $\alpha_1\beta\in B(11)$, $1\leq n\leq l-2$},
\\ 
(\phi) \otimes (1^l) & \text{if $\alpha_1\beta\in B(11)$, $n=l-1$},
\\ 
(1) \otimes (1^l) & \text{if $\alpha_1\beta\in B(11)$, $n=l$}, 
\\ 
(1) \otimes (\bar{1}) & \text{if $(\alpha_1,\beta)=(1,\bar{1})$, $n=1$},
\\ 
(1) \otimes (\phi) & \text{if $(\alpha_1,\beta)=(1,\bar{1})$, $n=2$}, 
\\ 
(1) \otimes (1^{n-2}) 
& \text{if $(\alpha_1,\beta)=(1,\bar{1})$, $3\leq n\leq l$}, 
\\
(1) \otimes (1^{n-1}0) & \text{for 
$\alpha_1\beta\in B(10)$, $\alpha_2\gamma\in B(11)$, $1\leq n\leq l$}, 
\\ 
(1) \otimes (1^{n-2}\bar{3}) & \text{for 
$\alpha_1\beta\in B(10)$, $\alpha_2\gamma\in B(12)$, $2\leq n\leq l$}, 
\\ 
(1) \otimes (1^{n-1}20) 
& \text{if $\alpha_1\beta\in B(12)$, $1\leq n\leq l-1$},
\\ 
(1) \otimes (1^{l-1}2) & \text{if $\alpha_1\beta\in B(12)$, $n=l$}, 
\\ 
(1)\otimes (\phi) & \text{if $b_1\neq (\phi)$, $b_2=(\phi)$, $l=1$}, 
\\ 
(\phi) \otimes (1^n) & \text{if $b_1\neq (\phi)$, $b_2=(\phi)$, 
$1\leq n\leq l-1$}, 
\\ 
(1) \otimes (1^{l-1}) & \text{if $b_1\neq (\phi)$, $b_2=(\phi)$, $n=l$}, 
\\ 
(\phi)\otimes (1) & \text{if $b_1=(\phi)$, $b_2\neq (\phi)$, $l=1$}, 
\\ 
(1)\otimes (1\bar{1}) & \text{if $b_1=(\phi)$, $b_2\neq (\phi)$, $l\neq 1$}, 
\\ 
(\phi)\otimes (\phi) & \text{if $b_1=(\phi)$, $b_2=(\phi)$}. 
\end{cases}
\end{align*}
Thanks to Lemma \ref{lem:comb.R for h.w.e.} and  Lemma \ref{lem:G2}, 
we have $\iota(b_1\otimes b_2)=b'_2\otimes b'_1$.  

Therefore the proof is finished.

\section{Soliton cellular automaton}\label{sec:SCA}

\subsection{States and time evolution}

Let $\{B_l\}_{l\geq 1}$ be the family of $U'_q(D_4^{(3)})$-crystals. 
Fix a sufficiently large positive integer $L$. 
Put $u_l=(l,0,0,0,0,0)\in B_l$. 
We define a set of paths $\mathscr{P}_L$ by 
\begin{equation}
\mathscr{P}_L=\left\{ \ 
p=b_1\otimes b_2\otimes \cdots \otimes b_L \in B_1^{\otimes L} 
;\ b_j=u_1 \ \text{if $j\gg 1$} \ \right\}. 
\end{equation}

By the iterating combinatorial $R$: 
$\text{Aff}(B_l)\otimes \text{Aff}(B_1)\simeq \text{Aff}(B_1)\otimes \text{Aff}(B_l)$, we have 
\begin{align*}
z^{0}b^{(0)}\otimes z^{0}b_1&\otimes z^{0}b_2 \otimes \cdots \otimes z^{0}b_L 
\\ \mapsto & \ 
z^{H_1}\tilde{b}_1 \otimes z^{-H_1}b^{(1)}\otimes z^{0}b_2 
\otimes  \cdots \otimes z^{0}b_L
\\ \mapsto & \ 
z^{H_1}\tilde{b}_1 \otimes z^{H_2}\tilde{b}_2 \otimes z^{-(H_1+H_2)}b^{(2)}\otimes z^{0}b_3\otimes \cdots \otimes z^{0}b_L 
\\ \mapsto & \ 
\cdots \qquad \cdots \qquad \cdots 
\\ \mapsto & \ \ 
z^{H_1}\tilde{b}_1 \otimes z^{H_2}\tilde{b}_2\otimes \cdots \otimes 
z^{H_L}\tilde{b}_L \otimes z^{E_l(b_1\otimes b_2\otimes \cdots \otimes b_L)}b^{(L)}, 
\end{align*}
\begin{equation*}
E_l(b_1\otimes b_2\otimes \cdots \otimes b_L)=-\sum_{j=1}^L H_j,
\quad H_j=H(b^{(j-1)}\otimes b_j). 
\end{equation*}
Here $b^{(L)}\in B_l$ is uniquely determined from 
$b^{(0)}\in B_l$ and $p\in \mathscr{P}_L$. 

We denote the iterating crystal isomorphism 
$B_l\otimes B_1\simeq B_1\otimes B_l$ by the following symbol. 
\[
\batten{b^{(0)}}{b_1}{b'_1}{}{a}
\hspace{-1mm}
\batten{b^{(1)}}{b_2}{b'_2}{}{a}
\hspace{-1mm}
\batten{b^{(2)}}{b_3}{b'_3}{}{a}
\hspace{-1mm}
\batten{b^{(3)}}{\cdots \ \cdots}{\cdots \ \cdots}{}{a}
\hspace{0.5cm}
\batten{b^{(L-1)}}{b_L}{b'_L}{b^{(L)}}{a}
\]

Set $b^{(0)}=u_l$. 
Assume $p=b_1\otimes b_2\otimes \cdots \otimes b_L$ is the element of 
$\mathscr{P}_L$. Then, 
\begin{align*}\label{eq:R matrix array}
b^{(L)}=u_l,\quad 
\tilde{b}_1\otimes \tilde{b}_2\otimes \cdots \otimes \tilde{b}_L 
\in \mathscr{P}_L, 
\quad 
E_l(p)<\infty. 
\end{align*}


Thus we define a map $T_l : \mathscr{P}_L\longrightarrow \mathscr{P}_L$ by 
\begin{equation}
T_l( b_1 \otimes b_2 \otimes \cdots \otimes b_L )=
\tilde{b}_1 \otimes \tilde{b}_2 \otimes \cdots \otimes \tilde{b}_L. 
\end{equation}
We call the dynamical system 
$\{T_l^t(p)\ \vert \ t\geq 0,\ l\geq 1, p\in \mathscr{P}_L \}$ 
the $D_4^{(3)}$-automaton.

\begin{example}\label{ex:time evolution}
Using the algorithm given in Theorem 4.4, we have 
\[
\batten{111}{\bar{2}}{1}{}{a}
\hspace{-2.5mm}
\batten{13}{0}{1}{}{a}
\hspace{-2.5mm}
\batten{1\bar{2}}{\bar{3}}{3}{}{a}
\hspace{-2.5mm}
\batten{2\bar{1}}{1}{\phi}{}{a}
\hspace{-1.5mm}
\batten{12\bar{1}}{1}{\bar{1}}{}{a}
\hspace{-1.5mm}
\batten{112}{1}{2}{}{a}
\hspace{-1.5mm}
\batten{111}{1}{1}{}{a}
\hspace{-1.5mm}
\batten{111}{1}{1}{\cdots}{a}
\]
Here the frames are omitted. Thus the time evolution is given by 
\[
T_3: \bar{2}\otimes 0\otimes \bar{3}\otimes 1\otimes 1\otimes 1 \otimes 1\otimes \cdots \ \mapsto \ 1\otimes 1\otimes 3\otimes \phi \otimes \bar{1}\otimes 2 \otimes 1\otimes \cdots.  
\]
\end{example}

\begin{lemma}
For a fixed element of $\mathscr{P}_L$, there exists an integer $l_0$ 
such that $T_l=T_{l_0}$ for any $l>l_0$. 
\end{lemma}

\begin{proposition}\label{prop:conservation laws}
(Conservation laws).\ 
For any element $p\in \mathscr{P}_L$, we have 
\begin{align}
T_lT_{l'}(p)=T_{l'}T_l(p), \quad 
E_l(T_{l'}(p))=E_l(p). 
\end{align}
\end{proposition}
\begin{proof}
Thanks to the Yang-Baxter equation, 
we can show this by the same argument as in Proposition 2.3 of [7]. 
\end{proof}

\subsection{Solitons}

In the $D_4^{(3)}$-automaton 
a state of the following form is called an $m$-soliton state of length 
$(l_1,l_2,\cdots ,l_m)$ 
\begin{equation}
...[l_1]......[l_2].....\cdots .....[l_m]..............
\end{equation}
Here $\ldots [l]\ldots$ denotes a local configuration such as 
\begin{equation}\label{eq:A_1-soliton}
\cdots \otimes (1)\otimes (1)
\otimes (3)^{\otimes x_2}\otimes (2)^{\otimes x_1}\otimes 
(1)\otimes (1)\otimes \cdots , \quad x_1+x_2=l
\end{equation}
sandwiched by sufficiently many $(1)$'s.

\begin{remark}
A one-soliton state in the $A_n^{(1)}$-automaton is characterized by the condition $E_1(p)=1$. 
We also adopt the definition of one-soliton state in $D_4^{(3)}$-automaton 
from it.  
\end{remark}

\begin{lemma}\label{lem:soliton step}
Let $p$ be a one soliton state of length $l$, then 
\begin{align}
& 
\text{The $k$-th conserved quantity of $p$ is given by $E_k(p)=\min(k,l)$.}
\\ & 
\text{$T_k(p)$ is obtained by the rightward shift by $\min(k,l)$ lattice steps.} 
\end{align}
\end{lemma}

Under the time evolution $T_r$
we identify the local state $``\ldots [l]\ldots "$ at time $t$ with 
\begin{equation}
\begin{array}{c}
z^{\gamma}(x_1,x_2);\ x_1\geq 0,\ x_2\geq 0,\ x_1+x_2=l, 
\\ {} \\ 
\gamma=\min(r,l)t+
\text{(the position of $``...[l]..."$ from the right end)}. 
\end{array}
\end{equation}
Here $\gamma$ denotes the phase of $``\ldots [l]\ldots "$. 
From Lemma \ref{lem:soliton step}, the phase 
$\gamma$ is invariant unless it interacts with other solitons.

Let $\tilde{B}_l$ ($l\in \mathbb{Z}_{\geq 1}$) be a 
$U'_q(A_{1}^{(1)})$-crystal, namely 
\begin{equation}
\tilde{B}_l=\{(x_1, x_2) \ \vert \ x_1\geq 0,\ x_2\geq 0,\ x_1+x_2=l\}. 
\end{equation}

\begin{proposition}
Define a map 
$\iota_l: \tilde{B}_l\sqcup \{0\}\rightarrow B_1^{\otimes l}\sqcup \{0\}$ by 
\begin{equation*}
\iota_l(b)=(3)^{\otimes x_2}\otimes (2)^{\otimes x_1} 
\quad \text{if $b\neq 0$} \ ;\ =0 \quad \text{otherwise}.
\end{equation*}
{} Then we have 
\begin{equation*}
\iota_l(\tilde{e}_1(b))=\tilde{e}_{2}(\iota_l(b)),\quad 
\iota_l(\tilde{f}_1(b))=\tilde{f}_{2}(\iota_l(b)) 
\quad \text{for any $b\in \tilde{B}_l\sqcup \{0\}$. }
\end{equation*}
\end{proposition}
{} Thus each local state $\ldots [l]\ldots$ is labeled by $U_q(A_1^{(1)})$-crystal $\text{Aff}(\tilde{B}_l)$. 


\begin{remark}
When $\mathfrak{g}_n$ is a non-exceptional algebra, a one-soliton state is labeled by the smaller algebra $U_q(\mathfrak{g}_{n-1})$. 
\end{remark}

\begin{example}
Using the algorithm of combinatorial $R$ given in Theorem 4.4, we have 
\[
\batten{111}{3}{1}{}{a}
\hspace{-1mm}
\batten{113}{3}{1}{}{a}
\hspace{-1mm}
\batten{133}{2}{1}{}{a}
\hspace{-1mm}
\batten{233}{1}{3}{}{a}
\hspace{-1mm}
\batten{123}{1}{3}{}{a}
\hspace{-1mm}
\batten{112}{1}{2}{}{a}
\hspace{-1mm}
\batten{111}{1}{1}{\cdots}{a}
\]
Thus the array ``$332$"  behaves as a soliton just like 
in the classical soliton theory. 
\end{example}

We denote a $m$-soliton state with length $(l_1,l_2,...,l_m)$ by 
\begin{equation}
z^{\gamma_1}b_1\otimes z^{\gamma_2}b_2\otimes \cdots 
\otimes z^{\gamma_m}b_m \in
\text{Aff}(\tilde{B}_{l_1})\otimes 
\text{Aff}(\tilde{B}_{l_2})\otimes \cdots \otimes \text{Aff}(\tilde{B}_{l_m}). 
\end{equation}

Assuming $l_1>l_2>\cdots >l_m$, 
we can expect that the state turns out to be 
\begin{equation}
................[l_m]........\cdots .....................[l_2]......[l_1]....
\end{equation}
after a sufficiently long time because the longer soliton moves faster 
under the time evolution $T_r$ $(r\gg 1)$.  
We describe such a scattering of solitons as follows.  
\begin{eqnarray*}
S_m: \text{Aff}(\tilde{B}_{l_1})\otimes \text{Aff}(\tilde{B}_{l_2})\otimes 
\cdots \otimes \text{Aff}(\tilde{B}_{l_m}) 
& \longrightarrow & 
\text{Aff}(\tilde{B}_{l_m})\otimes \cdots \otimes \text{Aff}(\tilde{B}_{l_2})
\otimes \text{Aff}(\tilde{B}_{l_1}) 
\\ 
z^{\gamma_1}b_1\otimes z^{\gamma_2}b_2 \otimes \cdots \otimes z^{\gamma_m}b_m
& \mapsto & 
z^{\gamma'_m}b'_m \otimes \cdots \otimes z^{\gamma'_2}b'_2 \otimes 
z^{\gamma'_1}b'_1
\end{eqnarray*}

The map $S_2$ signifies the two-body scattering of solitons. 
Now the theorem is 

\begin{theorem}\label{thm:two body scattering}
Assume $l_1>l_2$. Under the time evolution $T_r$ $(r>l_2)$ the map $S_2$ 
is described by the combinatorial $R$ matrix for $U_q(A_1^{(1)})$-crystals. 
Especially, the phase shift is given by 
\begin{equation}
\gamma'_2-\gamma_2=\gamma_1-\gamma'_1=2l_2+3\times H(b_1\otimes b_2). 
\end{equation} 
\end{theorem}

\begin{remark}
The combinatorial $R$ matrix in Theorem \ref{thm:two body scattering} has an extra term $2l_2$ in the power of $z$. 
However, the Yang-Baxter equation holds as it is. 
\end{remark}

\begin{corollary}\label{cor:multi-soliton scattering}
Scattering of solitons is factorized into two-body scattering. 
\end{corollary}
{} The map decomposes into the combinatorial $R$'s corresponding to pairwise transpositions of the components. It is independent of the order of the transpositions due to the Yang-Baxter equation. 
For example, we consider the $3$-body scattering of solitons with length 
$(l_1,l_2,l_3)$. 
Assume $l_1>l_2>l_3$. Under the time evolution $T_r$ $(r>l_1)$, 
the rule of 3-body scattering is given by 
\begin{equation}
S_3: 
z^{\gamma_1}b_1\otimes z^{\gamma_2}b_2\otimes z^{\gamma_3}b_3 
\mapsto z^{\gamma'_3}b'_3\otimes z^{\gamma'_2}b'_2\otimes z^{\gamma'_1}b'_1 
\end{equation}
The map $S_3$ is factorized into 
$(S_2\otimes 1)(1\otimes S_2)(S_2\otimes 1)$ or 
$(1\otimes S_2)(S_2\otimes 1)(1\otimes S_2)$. 
Since $S_2=R$ satisfies the Yang-Baxter equation, we have 
\begin{align*}
& 
(S_2\otimes 1)(1\otimes S_2)(S_2\otimes 1)
(z^{\gamma_1}b_1\otimes z^{\gamma_2}b_2\otimes z^{\gamma_3}b_3)
\\ & \qquad \qquad \qquad 
=(1\otimes S_2)(S_2\otimes 1)(1\otimes S_2)
(z^{\gamma_1}b_1\otimes z^{\gamma_2}b_2\otimes z^{\gamma_3}b_3). 
\end{align*}

\subsection{Examples}

Let us present examples of soliton scattering. 
Fix $L=50$.  
For a given state $p\in \mathscr{P}_L$ 
compute a series of calculations of $B_r\otimes B_1\simeq B_1\otimes B_r$ 
by the algorithm in Theorem \ref{thm:Combinatorial R}, and obtain 
the time evolution $T_r(p)$. 
For a given two-soliton state 
$z^{\gamma_1}(x_1,x_2)\otimes z^{\gamma_2}(y_1,y_2)\in 
\text{Aff}(\tilde{B}_{l_1})\otimes \text{Aff}(\tilde{B}_{l_2})$, 
put $H=2l_2-3\min(x_1,y_2)$. Let $\delta$ be the phase shift.

\begin{example}
The $r=4$ case. 
\[
\begin{array}{rl}
t=0:& \underline{2222}11111\underline{222}11111111111111111111111111111111111111
\\ 
t=1:& 1111\underline{2222}1111\underline{222}11111111111111111111111111111111111
\\ 
t=2:& 11111111\underline{2222}111\underline{222}11111111111111111111111111111111
\\ 
t=3:& 111111111111\underline{2221}11\underline{222}21111111111111111111111111111
\\ 
t=4:& 1111111111111112\underline{2211}1\underline{122}22111111111111111111111111
\end{array}
\]
\[
S_2=R: z^{46}(4,0)\otimes z^{38}(3,0)\mapsto 
z^{38+H}(3,0)\otimes z^{46-H}(4,0),
\quad \delta=H=6. 
\]
\end{example}

\begin{example}
In $r=3$ case. 
\[
\begin{array}{rl}
t=0:& \underline{222}1111\underline{33}11111111111111111111111111111111111111111
\\  
t=1:& 111\underline{222}111\underline{33}111111111111111111111111111111111111111
\\ 
t=2:& 111111\underline{222}11\underline{33}1111111111111111111111111111111111111
\\ 
t=3:& 111111111\underline{222}1\underline{33}11111111111111111111111111111111111
\\ 
t=4:& 111111111111\underline{222}\underline{33}111111111111111111111111111111111
\\ 
t=5:& 111111111111111\underline{2203}1111111111111111111111111111111
\\ 
t=6:& 111111111111111111\underline{2\underline{\bar{3}3}}11111111111111111111111111111
\\ 
t=7:& 111111111111111111111\underline{\underline{\bar{3}0}1}11111111111111111111111111
\\ 
t=8:& 11111111111111111111111\underline{\phi\bar{3}11}11111111111111111111111
\\ 
t=9:& 1111111111111111111111111\underline{1\bar{1}}\underline{\bar{2}21}11111111111111111111
\\ 
t=10:& 111111111111111111111111111\underline{1\phi}0\underline{211}11111111111111111
\\ 
t=11:& 11111111111111111111111111111\underline{11}\bar{3}3\underline{211}11111111111111
\\ 
t=12:& 1111111111111111111111111111111\underline{11}203\underline{211}11111111111
\\ 
t=13:& 111111111111111111111111111111111\underline{11}2233\underline{211}11111111
\\ 
t=14:& 11111111111111111111111111111111111\underline{11}22133\underline{211}11111
\\ 
t=15:& 1111111111111111111111111111111111111\underline{11}221133\underline{211}11
\end{array}
\]
\[
S_2=R: z^{47}(3,0)\otimes z^{41}(2,0)\mapsto 
z^{41+H}(2,0)\otimes z^{47-H}(1,2),
\quad \delta=H=-2. 
\]
\end{example}

\begin{example}\label{ex:3-body scattering}
The $r=3$ case. 
\[
\begin{array}{rl}
t=0:& 
33211132111311111111111111111111111111111111111111
\\ 
t=1:& 
11133211321131111111111111111111111111111111111111
\\ 
t=2:& 
11111133213213111111111111111111111111111111111111
\\ 
t=3:& 
11111111133232311111111111111111111111111111111111
\\ 
t=4:& 
11111111111133001111111111111111111111111111111111
\\ 
t=5:& 
1111111111111113\bar{1}311111111111111111111111111111111
\\ 
t=6:& 
111111111111111131\phi 0311111111111111111111111111111
\\ 
t=7:& 
111111111111111113111\bar{3}3311111111111111111111111111
\\ 
t=8:& 
11111111111111111131111203311111111111111111111111
\\ 
t=9:& 
11111111111111111113111112233311111111111111111111
\\ 
t=10:& 
11111111111111111111311111122133311111111111111111
\\ 
t=11:& 
11111111111111111111131111111221133311111111111111
\\ 
t=12:& 
11111111111111111111113111111112211133311111111111
\end{array}
\]
\[
S_3: z^{47}(1,2)\otimes z^{42}(1,1)\otimes z^{38}(0,1)\mapsto 
z^{37}(0,1)\otimes z^{41}(2,0)\otimes z^{47}(0,3)
\]
\begin{align*}
z^{47}(1,2)\otimes z^{42}(1,1)\otimes z^{38}(0,1)
\xrightarrow{R\otimes 1} & 
z^{42+1}(1,1)\otimes z^{47-1}(1,2)\otimes z^{38}(0,1)
\\ \xrightarrow{1\otimes R} & 
z^{43}(1,1)\otimes z^{38+(-1)}(1,0)\otimes z^{46-(-1)}(0,3)
\\ \xrightarrow{R\otimes 1} & 
z^{37+2}(0,1)\otimes z^{43-2}(2,0)\otimes z^{47}(0,3)
\end{align*}
\begin{align*}
z^{47}(1,2)\otimes z^{42}(1,1)\otimes z^{38}(0,1)
\xrightarrow{1\otimes R} & 
z^{47}(1,2)\otimes z^{38+(-1)}(1,0)\otimes z^{42-(-1)}(0,2)
\\ \xrightarrow{R\otimes 1} & 
z^{37+2}(0,1)\otimes z^{47-2}(2,1)\otimes z^{43}(0,2)
\\ \xrightarrow{1\otimes R} & 
z^{39}(0,1)\otimes z^{43+(-2)}(2,0)\otimes z^{45-(-2)}(0,3)
\end{align*}
\end{example}

\subsection{Operator $T_{\natural}$}

Let $B_{\natural}$ be the crystal given in section \ref{sec:Bnatural}. 
Put 
$u_{\natural}=
\begin{pmatrix}
1 \\ 2 
\end{pmatrix}
\in B_{\natural}$.

By the iterating crystal isomorphism 
$B_{\natural}\otimes B_1\rightarrow B_1\otimes B_{\natural}$ we have 
\begin{align*}
u_{\natural}\otimes b_1\otimes \cdots \otimes b_L 
\mapsto 
\tilde{b}_1 \otimes u_{\natural}^{(1)}\otimes b_2 \otimes \cdots \otimes b_L
\mapsto \cdots \mapsto  
\tilde{b}_1 \otimes \cdots \otimes \tilde{b}_L \otimes \tilde{u}_{\natural}. 
\end{align*}
Here $\tilde{u}_{\natural}\in B_{\natural}$ is uniquely determined from 
$p=b_1\otimes \cdots \otimes b_L\in \mathscr{P}_L$. 

Define a map 
$T_{\natural}: B_1^{\otimes L}\rightarrow B_1^{\otimes L}$ by 
\begin{equation} 
T_{\natural}(b_1\otimes b_2\otimes \cdots \otimes b_L)=
\tilde{b}_1\otimes \tilde{b}_2\otimes \cdots \otimes \tilde{b}_L. 
\end{equation}

Now we identify the two-soliton state 
\begin{align*}
p=(1)^{\otimes c}\otimes (3)^{\otimes x_2}\otimes (2)^{\otimes x_1}\otimes 
(1)^{\otimes d}\otimes (3)^{\otimes y_2}\otimes (2)^{\otimes y_1}\otimes (1)
\otimes \cdots \in \mathscr{P}_L,  
\\ 
(x_1+x_2=l,\  y_1+y_2=k,\ d\gg 1)
\end{align*}
with the tensor product of the element of $U_q(A_1^{(1)})$-crystal  
\begin{align}
z^{L-c-l}(x_1,x_2)\otimes z^{L-c-l-d-k}(y_1,y_2)
\in \text{Aff}(\tilde{B}_l)\otimes \text{Aff}(\tilde{B}_k). 
\end{align}

\begin{lemma}\label{lem:ironuki}
For the two-soliton state with highest weight 
$z^{\gamma_1}(l,0)\otimes z^{\gamma_2}(y_1,y_2)\in 
\text{Aff}(\tilde{B}_l)\otimes \text{Aff}(\tilde{B}_k)$  we have 
\[
T_{\natural}\left(
z^{\gamma_1}(l,0)\otimes z^{\gamma_2}(y_1,y_2)
\right)=
\begin{cases}
z^{\gamma_1}(l,0)\otimes z^{\gamma_2}(k,0) & \text{if $y_2=0$}, 
\\ 
z^{\gamma_1}(l,0)\otimes z^{\gamma_2-3}(y_1+1,y_2-1) 
& \text{if $y_2 \neq 0$}. 
\end{cases} 
\]
\end{lemma}
Instead of proving the lemma, we give examples below. 
(See also Appendix B.)
\begin{example}\label{ex:ironuki}
Omitting the frames, 
$T_{\natural}(z^{\gamma}(x_1,x_2))$ is described as follows. 
\end{example}
\begin{itemize}
\item
$T_{\natural}(z^{\gamma}(4,0))=z^{\gamma}(4,0)$. 
\[
\batten{{\begin{array}{c} 1 \\[-1mm] 2 \end{array}}}{1}{1}{}{a}
\hspace{-1.5mm}
\batten{{\begin{array}{c} 1 \\[-1mm] 2 \end{array}}}{1}{1}{}{a}
\hspace{-1.5mm}
\batten{{\begin{array}{c} 1 \\[-1mm] 2 \end{array}}}{1}{1}{}{a}
\hspace{-1.5mm}
\batten{{\begin{array}{c} 1 \\[-1mm] 2 \end{array}}}{2}{2}{}{a}
\hspace{-1.5mm}
\batten{{\begin{array}{c} 1 \\[-1mm] 2 \end{array}}}{2}{2}{}{a}
\hspace{-1.5mm}
\batten{{\begin{array}{c} 1 \\[-1mm] 2 \end{array}}}{2}{2}{}{a}
\hspace{-1.5mm}
\batten{{\begin{array}{c} 1 \\[-1mm] 2 \end{array}}}{1}{1}{}{a}
\hspace{-1.5mm}
\batten{{\begin{array}{c} 1 \\[-1mm] 2 \end{array}}}{1}{1}
{{\begin{array}{c} 1 \\[-1mm] 2 \end{array}}}{a}
\]

\item
$T_{\natural}(z^{\gamma}(3,1))=z^{\gamma}(4,0)$. 
\[
\batten{{\begin{array}{c} 1 \\[-1mm] 2 \end{array}}}{3}{1}{}{a}
\hspace{-1.5mm}
\batten{1_2}{2}{1}{}{a}
\hspace{-1.5mm}
\batten{2_1}{2}{1}{}{a}
\hspace{-1.5mm}
\batten{{\begin{array}{c} 2 \\[-1mm] \bar{3} \end{array}}}{2}{2}{}{a}
\hspace{-1.5mm}
\batten{{\begin{array}{c} 2 \\[-1mm] \bar{3} \end{array}}}{1}{2}{}{a}
\hspace{-1.5mm}
\batten{{\begin{array}{c} 2 \\[-1mm] 0 \end{array}}}{1}{2}{}{a}
\hspace{-1.5mm}
\batten{{\begin{array}{c} 2 \\[-1mm] 3 \end{array}}}{1}{2}{}{a}
\hspace{-1.5mm}
\batten{{\begin{array}{c} 2 \\[-1mm] 3 \end{array}}}{1}{1}
{{\begin{array}{c} 1 \\[-1mm] 3 \end{array}}}{a}
\]

\item
$T_{\natural}(z^{\gamma}(0,4))=z^{\gamma}(1,3)$. 
\[
\batten{{\begin{array}{c} 1 \\[-1mm] 2 \end{array}}}{3}{1}{}{a}
\hspace{-1.5mm}
\batten{1_2}{3}{1}{}{a}
\hspace{-1.5mm}
\batten{3_1}{3}{1}{}{a}
\hspace{-1.5mm}
\batten{{\begin{array}{c} 3 \\[-1mm] \bar{2} \end{array}}}{3}{3}{}{a}
\hspace{-1.5mm}
\batten{{\begin{array}{c} 3 \\[-1mm] \bar{2} \end{array}}}{1}{3}{}{a}
\hspace{-1.5mm}
\batten{{\begin{array}{c} 3 \\[-1mm] 0 \end{array}}}{1}{3}{}{a}
\hspace{-1.5mm}
\batten{{\begin{array}{c} 2 \\[-1mm] 3 \end{array}}}{1}{2}{}{a}
\hspace{-1.5mm}
\batten{{\begin{array}{c} 1 \\[-1mm] 3 \end{array}}}{1}{1}
{{\begin{array}{c} 1 \\[-1mm] 3 \end{array}}}{a}
\]
\end{itemize}

\begin{lemma}\label{lem:Tnatural R=R Tnatural}
Assume that $l>k$. For the two-soliton state with highest weight 
$z^{\gamma_1}b_1\otimes z^{\gamma_2}b_2\in 
\text{Aff}(\tilde{B}_l)\otimes \text{Aff}(\tilde{B}_k)$ we have 
\[
T_{\natural}(R(z^{\gamma_1}b_1\otimes z^{\gamma_2}b_2))=
R(T_{\natural}(z^{\gamma_1}b_1\otimes z^{\gamma_2}b_2)). 
\]
\end{lemma}
\begin{proof} Put 
$z^{\gamma_1}b_1\otimes z^{\gamma_2}b_2=
z^{\gamma_1}(l,0)\otimes z^{\gamma_2}(y_1,y_2)\in 
\text{Aff}(\tilde{B}_l)\otimes \text{Aff}(\tilde{B}_k)$. 
The $y_2=0$ case is obvious. Assume $y_2\neq 0$, then we have  
\begin{align*}
R(z^{\gamma_1}b_1\otimes z^{\gamma_2}b_2)=&
z^{\gamma_2+(2k-3y_2)}(k,0)\otimes z^{\gamma_1-(2k-3y_2)}(l-y_2,y_2), 
\\ 
T_{\natural}(R(z^{\gamma_1}b_1\otimes z^{\gamma_2}b_2))=&
z^{\gamma_2+2k-3y_2}(k,0)\otimes z^{\gamma_1-2k+3y_2-3}(l-y_2+1,y_2-1). 
\\ 
T_{\natural}(z^{\gamma_1}b_1\otimes z^{\gamma_2}b_2)=&
z^{\gamma_1}(l,0)\otimes z^{\gamma_2-3}(y_1+1,y_2-1), 
\\ 
R(T_{\natural}(z^{\gamma_1}b_1\otimes z^{\gamma_2}b_2))=&
z^{(\gamma_2-3)+(2k-3(y_2-1))}(k,0)\otimes 
z^{\gamma_1-(2k-3(y_2-1))}(l-y_2+1,y_2-1).  
\end{align*}
Hence we have $T_{\natural}R=RT_{\natural}$. 
\end{proof}

\begin{lemma}\label{lem:Tnatural Tl=Tl Tnatural}
$T_{\natural}(T_l(p))=T_l(T_{\natural}(p))$ 
for any $l>0$, $p\in \mathscr{P}_L$
\end{lemma}
\begin{proof} 
It follows from the Yang-Baxter equation on 
$\text{Aff}(B_{\natural})\otimes \text{Aff}(B_l)\otimes \text{Aff}(B_{l'})$. 
\end{proof}

\begin{lemma}\label{lem:identiti}
For any $l,k\in \mathbb{Z}_{>0}$, 
let $z^{\gamma_1}b_1\otimes z^{\gamma_2}b_2$ and 
$z^{\gamma'_1}b'_1\otimes z^{\gamma'_2}b'_2$ 
be the highest weight elements in 
$\text{Aff}(\tilde{B}_l)\otimes \text{Aff}(\tilde{B}_k)$. 
Suppose that they satisfy the conditions 
$T_{\natural}(z^{\gamma_1}b_1\otimes z^{\gamma_2}b_2)=
T_{\natural}(z^{\gamma'_1}b'_1\otimes z^{\gamma'_2}b'_2)$ and  
$\text{wt}(b_1\otimes b_2)=\text{wt}(b'_1\otimes b'_2)$. Then we have 
$z^{\gamma_1}b_1\otimes z^{\gamma_2}b_2=
z^{\gamma'_1}b'_1\otimes z^{\gamma'_2}b'_2$.  
\end{lemma}
\begin{proof} Put 
$z^{\gamma_1}b_1\otimes z^{\gamma_2}b_2=
z^{\gamma_1}(l,0)\otimes z^{\gamma_2}(k-s,s)$ and  
$z^{\gamma'_1}b'_1\otimes z^{\gamma'_2}b'_2=
z^{\gamma'_1}(l',0)\otimes z^{\gamma_2}(k'-s',s')$. 
By weight condition to show is in the cases 
(i) $l+k=l'+k'$, $(s,s')=(0,0)$ or (ii) $l+k=l'+k'$, $(s,s')\neq (0,0)$.  
The case (i) is obvious. In (ii) case. 
By the condition 
$T_{\natural}(z^{\gamma_1}b_1\otimes z^{\gamma_2}b_2)=
T_{\natural}(z^{\gamma'_1}b'_1\otimes z^{\gamma'_2}b'_2)$ 
we have 
\[
z^{\gamma_1}(l,0)\otimes z^{\gamma_2-3}(k-s+1,s-1)=
z^{\gamma'_1}(l',0)\otimes z^{\gamma'_2-3}(k'-s'+1,s'-1). 
\]
Hence we have 
$(\gamma_1,\gamma_2,l,k,s)=(\gamma'_1,\gamma'_2,l',k',s')$. 
\end{proof}

\subsection{Proof of Theorem \ref{thm:two body scattering}}

For any two-soliton state 
$z^{\gamma_1}(x_1,x_2)\otimes z^{\gamma_2}(y_1,y_2)\in 
\text{Aff}(\tilde{B}_{l_1})\otimes \text{Aff}(\tilde{B}_{l_2})$, 
$T_r$ $(r\geq 1)$ commutes with $\tilde{e}_2, \tilde{f}_2$. 
Thus it is enough to check the rule for the highest weight elements 
$z^{\gamma_1}(l,0)\otimes z^{\gamma_2}(y_1,y_2)\in 
\text{Aff}(\tilde{B}_{l_1})\otimes \text{Aff}(\tilde{B}_{l_2})$. 
We will show the statement by the induction on $y_2$. 
Because of Proposition \ref{prop:conservation laws}, we have 
$T_r^bT_{l_2+1}^a=T_{l_2+1}^aT_r^b$ $(r>l_2)$. 
If $a$ and $b$ are sufficiently large, 
we can reduce the observation of the scattering by $T_r$ 
in the right hand side to that by $T_{l_2+1}$ in the left hand side. 
Hence the statement is reduced to $T_{l_2+1}$. 

The $y_2=0$ case. 
By direct calculation we obtain 
\begin{equation}
S_2: z^{\gamma_1}(l,0)\otimes z^{\gamma_2}(k,0)\mapsto 
z^{\gamma_2+2k}(k,0)\otimes z^{\gamma_1-2k}(l,0). 
\end{equation}

The $y_2>0$ case. 
For the two-soliton state 
\begin{align*}
p=&(1)^{\otimes (L-l_1-\gamma_1)}\otimes (2)^{\otimes l_1}\otimes 
(1)^{\otimes (\gamma_1-\gamma_2-l_2)}\otimes (3)^{\otimes y_2}\otimes (2)^{\otimes y_1}
\otimes (1)^{\otimes \gamma_2}\in \mathscr{P}_L
\\=& 
z^{\gamma_1}(l_1,0)\otimes z^{\gamma_2}(y_1,y_2)
\in \text{Aff}(\tilde{B}_{l_1})\otimes \text{Aff}(\tilde{B}_{l_2}), 
\end{align*}
the statement to show is the following: 
\begin{equation}\label{eq:final}
T_r^t(p)=R(p) \quad \text{for $r>l_2$, $t\gg 1$}. 
\end{equation}
{} Thanks to Lemma \ref{lem:Tnatural R=R Tnatural}  
and Lemma \ref{lem:Tnatural Tl=Tl Tnatural}, we have 
$(RT_{\natural})(p)=(T_{\natural}R)(p)$ and 
$(T_{\natural}T_r^t)(p)=(T_r^tT_{\natural})(p)$ respectively. 
By the induction hypothesis 
we have $(RT_{\natural})(p)=(T_r^tT_{\natural})(p)$. 
Hence we have $T_{\natural}(T_r^t(p))=T_{\natural}(R(p))$. 
Since both $T_r^t(p)$ and $R(p)$ have the same weight as $p$, 
we can apply Lemma \ref{lem:identiti}. Thus we have $(\ref{eq:final})$. 

Therefore the proof is finished.

\appendix

\section{Proof of Proposition \ref{lem:comb.R for h.w.e.}}\label{app:A}

Because of the weight, we have 
\begin{equation}
\iota :\ (l,0,0,0,0,0)\otimes (1,0,0,0,0,0)\mapsto 
(1,0,0,0,0,0)\otimes (l,0,0,0,0,0)
\end{equation}
under the crystal isomorphism 
$\iota: B_l\otimes B_1 \rightarrow  B_1\otimes B_l$. 

Set $u_n=(n,0,0,0,0,0)\in B_l$ for $0\leq n\leq l$. 
For any $b_1\otimes b_2\in B_l\otimes B_1$ there exist 
$\tilde{g}_{i_1},\ldots ,\tilde{g}_{i_m}
\in \{\tilde{e}_i, \tilde{f}_i \ \vert i=0,1,2\}$ such that 
$b_1\otimes b_2=\tilde{g}_{i_1}\cdots \tilde{g}_{i_m}(u_l\otimes u_1)$. 
Set $\tilde{b}_2\otimes \tilde{b}_1=\iota(b_1\otimes b_2)$. 
Then $\tilde{b}_2\otimes \tilde{b}_1=\tilde{g}_{i_1}\ldots\tilde{g}_{i_m}(u_1\otimes u_l)$. 
Denote it by $b_1\otimes b_2\simeq \tilde{b}_2\otimes \tilde{b}_1 ; \tilde{g}_{i_1}\cdots \tilde{g}_{i_m}$.

By using the following data: 
\begin{align}
\tilde{f}_0^p(0,0,0,0,0,l)=&
\begin{cases}
(0,0,0,0,0,l-p) & \text{if $0\leq p\leq l$}, 
\\ 
(p-l,0,0,0,0,0) & \text{if $l+1\leq p\leq 2l$}. 
\end{cases}
\\ 
\tilde{f}_0^p(0,0,1,1,0,l-1)=&
\begin{cases}
(0,0,1,1,0,l-1-p) & \text{if $0\leq p\leq l-1$}, 
\\ 
(p-l+1,0,1,1,0,0) & \text{if $l\leq p\leq 2l-2$}. 
\end{cases}
\\ 
\tilde{f}_0^p(0,0,0,2,0,l-1)=&
\begin{cases}
(0,0,0,2,0,l-1-p) & \text{if $0\leq p\leq l-1$}, 
\\ 
(p-l,1,0,0,0,0) & \text{if $l\leq p\leq 2l-1$}. 
\end{cases}
\\ 
\tilde{f}_0^p(0,1,0,0,0,l-1)=&
\begin{cases}
(0,1,0,0,0,l-1-p) & \text{if $0\leq p\leq l-2$}, 
\\ 
(p-l+1,1,1,1,0,0) & \text{if $l-1\leq p\leq 2l-3$}. 
\end{cases}
\\ 
\tilde{f}_0^p(1,0,0,0,0,l-1)=&
\begin{cases}
(1,0,0,0,0,l-1-p) & \text{if $0\leq p\leq l-2$}, 
\\ 
(p-l+3,0,0,0,0,1) & \text{if $l-1\leq p\leq 2l-4$}. 
\end{cases}
\end{align}
we obtain the crystal isomorphism 
$b_1\otimes b_2\simeq \tilde{b}_2\otimes \tilde{b}_1 ; \tilde{g}_{i_1}\cdots \tilde{g}_{i_m}$.  

The $l=1$ case. 
\[
\begin{array}{l}
(1)\otimes (\phi)\simeq (1)\otimes (\bar{1});\ \tilde{e}_0, 
\\ 
(1)\otimes (\bar{1})\simeq (1)\otimes (\bar{1});\ \tilde{e}_0^2, 
\\ 
(\phi)\otimes (1)\simeq (\phi)\otimes (1);\ E^1\tilde{e}_0^3, 
\end{array}
\begin{array}{l}
(\phi)\otimes (\phi)\simeq (\phi)\otimes (\phi);\
\tilde{e}_0E^1\tilde{e}_0^3, 
\\ 
(1)\otimes (0)\simeq (1)\otimes (0);\ 
\tilde{f}_0\tilde{f}_1\tilde{f}_2\tilde{f}_1E^1\tilde{e}_0^3,
\\ 
(1)\otimes (2)\simeq (1)\otimes (2);\
 \tilde{f}_0\tilde{f}_1\tilde{f}_0\tilde{f}_1\tilde{f}_2\tilde{f}_1E^1\tilde{e}_0^3. 
\end{array}
\]

The $l=2$ case. 
\[
\begin{array}{l}
(11)\otimes (\phi)\simeq (1)\otimes (1);\ \tilde{e}_0, 
\\ 
(11)\otimes (\bar{1})\simeq (1)\otimes (\phi);\ \tilde{e}_0^2, 
\\ 
(1)\otimes (\bar{1})\simeq (1)\otimes (\bar{1});\ \tilde{e}_0^3, 
\\ 
(1)\otimes (1)\simeq (\phi)\otimes (11);\ E^2\tilde{e}_0^5, 
\\ 
(1)\otimes (\phi)\simeq (\phi)\otimes (1);\ \tilde{e}_0E^2\tilde{e}_0^5, 
\\ 
(\phi)\otimes (\phi)\simeq (\phi)\otimes (\phi);\ 
\tilde{e}_0^2E^2\tilde{e}_0^5, 
\end{array}
\begin{array}{l}
(\phi)\otimes (1)\simeq (1)\otimes (1\bar{1});\ E^1\tilde{e}_0^4, 
\\ 
(1)\otimes (2)\simeq (1)\otimes (20);\ 
\tilde{f}_0\tilde{f}_1E^1\tilde{e}_0^4 
\\ 
(1)\otimes (0)\simeq (1)\otimes (0);\ 
\tilde{f}_0\tilde{f}_1\tilde{f}_2\tilde{f}_1E^1\tilde{e}_0^4
\\
(11)\otimes (0)\simeq (1)\otimes (10);\ 
\tilde{f}_0^2\tilde{f}_1\tilde{f}_2\tilde{f}_1E^1\tilde{e}_0^4, 
\\
(11)\otimes (\bar{3})\simeq (1)\otimes (2);\ \tilde{f}_0\tilde{f}_1\tilde{f}_0\tilde{f}_1\tilde{f}_2\tilde{f}_1E^1\tilde{e}_0^4, 
\\ 
(11)\otimes (2)\simeq (1)\otimes (12);\ \tilde{f}_0^2\tilde{f}_1\tilde{f}_0\tilde{f}_1\tilde{f}_2\tilde{f}_1E^1\tilde{e}_0^4. 
\end{array}
\]

The $l\geq 3$ case. 
\[
\begin{array}{l}
(1^l)\otimes (\phi)\simeq (1)\otimes (1^{l-1});\ \tilde{e}_0
\\  
(1^{l-n})\otimes (\bar{1})\simeq (1)\otimes (1^{l-n-2});\ 
\tilde{e}_0^{n+2} \ \text{for} \ 0\leq n\leq l-2, 
\\ 
(1)\otimes (\bar{1})\simeq (1)\otimes (\bar{1});\ \tilde{e}_0^{l+1}, 
\\   
(\phi)\otimes (1)\simeq (1)\otimes (1\bar{1});\ E^1\tilde{e}_0^{l+2}, 
\\  
(1^n)\otimes (2)\simeq (1)\otimes (1^{n-1}20);\ 
\tilde{f}_0^n\tilde{f}_1E^1\tilde{e}_0^{l+2} \ \text{for}\ 1\leq n\leq l-1,
\\ 
(1^n)\otimes (0)\simeq (1)\otimes (1^{n-1}0);\ 
\tilde{f}_0^{n-1}\tilde{e}_1\tilde{f}_0\tilde{f}_1^2\tilde{f}_2\tilde{f}_1E^1\tilde{e}_0^{l+2} \ \text{for}\ 1\leq n\leq l, 
\\ 
(1^n)\otimes (\bar{3})\simeq (1)\otimes (1^{n-2}2);\ 
\tilde{f}_0^n\tilde{f}_1^2\tilde{f}_2\tilde{f}_1E^1\tilde{e}_0^{l+2}
\ \text{for}\ 2\leq n\leq l-1, 
\\ 
(1^l)\otimes (2)\simeq (1)\otimes (1^{l-1}2);\ 
\tilde{f}_0^{l+1}\tilde{f}_1^2\tilde{f}_2\tilde{f}_1E^1\tilde{e}_0^{l+2} 
\\ 
(1^n)\otimes (1)\simeq (1)\otimes (1^{n+1}\bar{1});\ 
E^{n+1}\tilde{e}_0^{l+n+2}\ \text{for}\ 1\leq n\leq l-2, 
\\ 
(1^{l-1})\otimes (1)\simeq (\phi)\otimes (1^l);\ E^l\tilde{e}_0^{2l+1},  
\\ 
(1^{l-n})\otimes (\phi)\simeq (\phi)\otimes (1^{l-n});\ 
\tilde{e}_0^nE^l\tilde{e}_0^{2l+1} \ \text{for}\ 1\leq n\leq l.  
\end{array}
\]
Here 
$E^n=\tilde{e}_1^n\tilde{e}_2^n\tilde{e}_1^{2n}\tilde{e}_2^n\tilde{e}_1^n$.

\section{The table of 
$B_{\natural}\otimes B_1\simeq B_1\otimes B_{\natural}$}\label{app:B}

The following data is all images of the crystal isomorphism 
$B_{\natural}\otimes B_1\rightarrow B_1\otimes B_{\natural}$ 
that is needed for Lemma \ref{lem:ironuki}. 
\[
\begin{array}{lll}
{}^t(12) \otimes (1)\simeq (1)\otimes {}^t(12) 
& 
{}^t(12) \otimes (2)\simeq (2)\otimes {}^t(12)  
& 
{}^t(12) \otimes (3)\simeq (1)\otimes  (1)_2
\\ 
(1)_2 \otimes (1)\simeq (1)\otimes (1)_1 
& 
(1)_1\otimes (1)\simeq (1)\otimes {}^t(23) 
& 
{}^t(23) \otimes (1)\simeq (2)\otimes  {}^t(13)
\\ 
{}^t(13)\otimes (1)\simeq (1)\otimes  {}^t(13) 
& 
(1)_2\otimes (2)\simeq (1)\otimes (2)_1, 
& 
(2)_1\otimes (1)\simeq (1)\otimes  {}^t(20) 
\\ 
{}^t(20) \otimes (1)\simeq (2)\otimes  {}^t(23) 
& 
(2)_1\otimes (2)\simeq (1)\otimes {}^t(2\bar{3}) 
& 
{}^t(2\bar{3}) \otimes (2)\simeq (2)\otimes  {}^t(2\bar{3}) 
\\ 
{}^t(2\bar{3})\otimes (1)\simeq (2)\otimes {}^t(20)  
& 
(1)_2\otimes (3)\simeq (1)\otimes (3)_1
& 
(3)_1\otimes (1)\simeq (1)\otimes {}^t(30)  
\\ 
{}^t(30) \otimes (1)\simeq (3)\otimes  {}^t(23) 
& 
(3)_1\otimes (2)\simeq (1)\otimes  {}^t(3\bar{3})  
& 
{}^t(3\bar{3}) \otimes (1)\simeq (3)\otimes  {}^t(20) 
\\ 
{}^t(3\bar{3}) \otimes (2)\simeq (3)\otimes {}^t(2\bar{3}) 
& 
(3)_1\otimes (3)\simeq (1)\otimes  {}^t(3\bar{2})  
& 
{}^t(3\bar{2}) \otimes (1)\simeq (3)\otimes {}^t(30)  
\\ 
{}^t(3\bar{2})  \otimes (3)\simeq (3)\otimes {}^t(3\bar{2}) 
& 
{}^t(3\bar{2})  \otimes (2)\simeq (3)\otimes {}^t(3\bar{3}). & 
\end{array}
\]
where ${}^t(\alpha\beta)$ denotes 
$\begin{pmatrix} \alpha \\ \beta \end{pmatrix}$. 


\begin{thebibliography}{KMN2}

\bibitem{1}
R.~J.~Baxter, 
\textit{Exactly solved models in statical mechanics}, 
Academic Press, London (1982). 

\bibitem{2}
K.~Fukuda, M.~Okado and Y.~Yamada, 
\textit{Energy functions in box ball systems}, 
Int. J. Mod. Phys. A {\bf 15} (2000), 1379-1392.

\bibitem{3}
G.~Hatayama, K.~Hikami, R.~Inoue, A.~Kuniba, T.~Takagi and T.~Tokihiro,
\textit{The $A_M ^{(1)}$ Automata related to crystals of symmetric tensors}, 
J. Math. Phys. 42 (2001) 274-308.

\bibitem{4}
G.~Hatayama, A.~Kuniba, M.~Okado, and T.~Takagi, 
\textit{Combinatorial $R$ matrices for a family of crystals:\ 
$C_n^{(1)}$ and $A_{2n-1}^{(2)}$ cases}, 
in ``Physical Combinatorics", Prog. in Math. M.~Kashiwara and T.~Miwa ed. 
Birkhauser (2000) 105-139. 

\bibitem{5}
G.~Hatayama, A.~Kuniba, M.~Okado, and T.~Takagi, 
\textit{Combinatorial $R$ matrices for a family of crystals:\ 
$B_n^{(1)}$, $D_n^{(1)}$, $A_{2n^{(2)}}$ and $D_{n+1}^{(2)}$ cases}, 
J. Algebra {\bf 247} (2002) 577-615. 


\bibitem{6} G.~Hatayama, A.~Kuniba, M.~Okado, T.~Takagi and Z.~Tsuboi,
\textit{Paths, crystals and fermionic formulae}, MathPhys Odyssey 2001,
205--272, Prog.\ Math.\ Phys.\ {\bf 23}, Birkh\"auser Boston, Boston, MA, 2002.



\bibitem{7}
G.~Hatayama, A.~Kuniba, M.~Okado, T.~Takagi and Y.~Yamada, 
\textit{Scattering rules in soliton cellular automata associated with crystal bases}, 
Contemp. Math. {\bf 297} (2002), 151-182.


\bibitem{8}
G.~Hatayama, A.~Kuniba and T.~Takagi,
\textit{Soliton cellular automata associated with crystal bases}, 
Nucl. Phys. B577[PM](2000) 619-645.


\bibitem{9} M.~Kashiwara,
\textit{On crystal bases of the $q$-analogue of universal enveloping algebras},
Duke Math.\ J.\ {\bf 63}, 465-516 (1991).


\bibitem{10}
V. G. Kac, 
\textit{``Infinite Dimensional Lie Algebras,"}
3rd ed., Cambridge Univ. Press, Cambridge, UK, 1990.



\bibitem{11} S-J,~Kang, M.~Kahiwara, K.C.~Misra,
\textit{Crystal bases of Verma Modules for Quantum Affine Lie Algebras}, 
Compositio Math.\ {\bf 92} (1994), 299--325


\bibitem{12} S-J.Kang, K.C.Misra, 
\textit{Crystal Bases and Tensor Product Decompositions of $U_q(G_2)$-Modules},  J. Algebra {\bf 163}, (1994), 675--691


\bibitem{13} M.~Kashiwara and T.Nakashima,
\textit{Crystal graphs for representations of q-analogue of classical Lie 
algebras}, J.\ Algebra {\bf 163} (1994), 295--345.


\bibitem{13}
S-J.~Kang, M.~Kashiwara, K.~C.~Misra, 
T.~Miwa, T.~Nakashima and A.~Nakayashiki,
\textit{Affine crystals and vertex models},
Int.\ J.\ Mod.\ Phys.\ A {\bf 7} (suppl. 1A), 449-484 (1992).


\bibitem{14}
S-J.~Kang, M.~Kashiwara, K.~C.~Misra,
T.~Miwa, T.~Nakashima and A.~Nakayashiki,
\textit{Perfect crystals of quantum affine Lie algebras},
Duke Math.\ J.\ {\bf 68} (1992) 499--607.



\bibitem{15}
M.~Kashiwara, K.~C.~Misra, M.~Okado, D.~Yamada, 
\textit{Perfect crystals for $U_q(D_4^{(3)})$}, 
arXiv:math.QA/0610873v1. 



\bibitem{16} C.~Lecouvey, 
Combinatorics of crystal graphs 
for the root systems of types $A_{n}$, $B_{n}$, $C_{n}$, $D_{n}$ and 
$G_{2}$ \`{a} MSJ memoirs vol 17 (2007). 



\bibitem{17} A.~Nakayashiki and Y.~Yamada, 
\textit{Kostka polynomials and energy functions in solvable lattice models},  
Selecta Mathematica, New Ser.\ {\bf 3} (1997) 547-345. 



\bibitem{18} D.~Takahashi,
\textit{On some soliton systems defined by using
boxes and balls}, Proceedings of
the International Symposium on Nonlinear Theory and
Its Applications (NOLTA '93),
(1993) 555--558.


\bibitem{19} T.~Tokihiro, A.~Nagai and J.~Satsuma,
\textit{Proof of solitonical nature of box and ball systems by
means of inverse ultra-discretization},
inverse Probl. 15(1999) 1639-1662.


\bibitem{19} D.~Takahashi and J.~Satsuma,
\textit{A soliton cellular automaton},
J. Phys. Soc. Jpn. {\bf 59} (1990) 3514--3519.


\bibitem{20} T.~Tokihiro, D.~Takahashi, J.~Matsukidaira 
and J.~Satsuma,
\textit{From soliton equations to integrable cellular automata through
a limiting procedure},
Phys. Rev. Lett. {\bf 76} (1996) 3247--3250.


\bibitem{21} S.~Yamane,
\textit{Perfect crystals of $U_q(G_2^{(1)})$}, 
J.\ Algebra {\bf 210} (1998), 440--486.


\bibitem{22} D.~Yamada,
\textit{Box ball system associated with antisymmetric tensor crystals}, 
J. Phys. A 37 (2004), no. 42, 9975--9987. 


\end{thebibliography}
\end{document}